\documentclass[a4paper,reqno]{amsart}
\usepackage{amsmath}
\usepackage{amsfonts}
\usepackage{amssymb}
\usepackage{bm}

\newcommand\eps{\varepsilon}
\newcommand\R{{\mathbf{R}}}
\newcommand\C{{\mathbf{C}}}
\newcommand\Z{{\mathbf{Z}}}

\newcommand\loc{{\operatorname{loc}}}

\newcommand\rad{{\operatorname{rad}}}

\theoremstyle{plain}
  \newtheorem{theorem}[subsection]{Theorem}
  \newtheorem{conjecture}[subsection]{Conjecture}
  \newtheorem{proposition}[subsection]{Proposition}
  \newtheorem{lemma}[subsection]{Lemma}
  \newtheorem{corollary}[subsection]{Corollary}

\theoremstyle{remark}
  \newtheorem{remark}[subsection]{Remark}
  
  \newtheorem{example}[subsection]{Example}

\theoremstyle{definition}
  \newtheorem{definition}[subsection]{Definition}

\include{psfig}

\begin{document}
\title[Minimal-mass blowup solutions]{Minimal-mass blowup solutions of the mass-critical NLS}
\author{Terence Tao}
\address{Department of Mathematics, UCLA, Los Angeles CA 90095-1555}
\email{tao@@math.ucla.edu}

\author{Monica Visan}
\address{Institute for Advanced Study}
\email{mvisan@@ias.edu}

\author{Xiaoyi Zhang}
\address{Academy of Mathematics and System Sciences, Chinese Academy of Sciences}
\email{zh.xiaoyi@@gmail.com}


\vspace{-0.3in}
\begin{abstract}
We consider the minimal mass $m_0$ required for solutions to the mass-critical nonlinear Schr\"odinger (NLS) equation $iu_t + \Delta u = \mu |u|^{4/d} u$
to blow up.  If $m_0$ is finite, we show that there exists a minimal-mass solution blowing up (in the sense of an infinite spacetime norm) in both time
directions, whose orbit in $L^2_x(\R^d)$ is compact after quotienting out by the symmetries of the equation.  A similar result is obtained for
spherically symmetric solutions.  Similar results were previously obtained by Keraani, \cite{keraani}, in dimensions $1$, $2$ and Begout and Vargas,
\cite{begout}, in dimensions $d\geq 3$ for the mass-critical NLS and by Kenig and Merle, \cite{merlekenig}, in the energy-critical case.
In a subsequent paper we shall use this compactness result to establish global existence and scattering in $L^2_x(\R^d)$ for the defocusing NLS
in three and higher dimensions with spherically symmetric data.
\end{abstract}

\maketitle

\section{Introduction}

\subsection{The mass-critical nonlinear Schr\"odinger equation}

Fix a dimension $d \geq 1$ and a sign $\mu = \pm 1$ ($\mu=+1$ is \emph{defocusing}, $\mu=-1$ is \emph{focusing}).  We shall consider strong
$L^2_x(\R^d)$ solutions to the mass-critical (or \emph{pseudoconformal}) nonlinear Schr\"odinger (NLS) equation
\begin{equation}\label{nls}
i u_t + \Delta u = F(u)
\end{equation}
where $F(z) := \mu |z|^{4/d} z$ is the nonlinearity.
More precisely, we say that a function $u: I \times \R^d \to \C$ on a time interval $I \subset \R$ (possibly half-infinite or infinite) is a
\emph{strong $L^2_x(\R^d)$ solution} (or \emph{solution} for short) to \eqref{nls} if it lies in the
class\footnote{We use $C^0_t L^2_x(I \times \R^d)$ to denote the Banach space of spacetime functions $u: I \times \R^d \to \C$
such that the curve $t \mapsto u(t)$ is bounded and continuous in $L^2_x$, with the uniform norm $\sup_{t \in I} \|u(t)\|_{L^2_x(\R^d)}$,
whereas $C^0_{t,\loc} L^2_x(I \times \R^d)$ is the Frechet space of functions where the curve is merely continuous, but not necessarily bounded.
Similarly for $L^{2(d+2)/d}_{t,x}(I \times \R^d)$ and $L^{2(d+2)/d}_{t,\loc} L^{2(d+2)/d}_x(I \times \R^d)$.}
$C^0_{t,\loc} L^2_x(I \times \R^d) \cap L^{2(d+2)/d}_{t,\loc} L^{2(d+2)/d}_x(I \times \R^d)$, and we have the Duhamel formula
$$
u(t_1) = e^{i(t_1-t_0)\Delta} u(t_0) - i \int_{t_0}^{t_1} e^{i(t_1-t)\Delta} F(u(t))\ dt
$$
for all $t_0, t_1 \in I$.   Here, $e^{it\Delta}$ is the propagator for the free Schr\"odinger equation, defined via the Fourier transform
$$
\hat f(\xi) := \int_{\R^d} e^{-ix \cdot \xi} f(x)\ dx
$$
by
$$
\widehat{e^{it\Delta} f}(\xi) = e^{-it|\xi|^2} \hat f(\xi).
$$
We say that the solution has \emph{maximal lifespan} if the interval $I$ cannot be extended to any strictly larger interval.
The condition $u \in L^{2(d+2)/d}_{t,\loc} L^{2(d+2)/d}_x$ is a natural one arising from the Strichartz perturbation theory; for instance,
it is currently necessary in order to ensure uniqueness of (maximal-lifespan) solutions.  Solutions to \eqref{nls} in this class have been
intensively studied, see e.g. \cite{begout}, \cite{bourg.2d}, \cite{ck}, \cite{caz}, \cite{cwI}, \cite{keraani}, \cite{merle}, \cite{merle2},
\cite{merle3}, \cite{mr}, \cite{mr2}, \cite{mr3}, \cite{mt}, \cite{mv}, \cite{tsutsumi}.

We will measure the size of such solutions in two different ways.  Firstly, we define the \emph{mass} $M(f)$ of a function $f \in L^2_x(\R^d)$ by
$$
M(f) := \|f\|_{L^2_x(\R^d)}^2 =\int_{\R^d} |f(x)|^2\ dx.
$$
We shall shortly see that if $u$ is a solution, then $M(u(t))$ is independent of $t$, and so we may meaningfully discuss the mass $M(u)$ of a solution.
Secondly, we define the \emph{scattering size} $S(u) = S_I(u) \in [0,+\infty]$ of a function $u: I \times\R^d \to \C$ (which, in particular, could be
a solution) by
$$
S(u) := \| u \|_{L^{2(d+2)/d}_{t,x}(I \times \R^d)}^{2(d+2)/d} = \int_I \int_{\R^d} |u(t,x)|^{2(d+2)/d}\ dx dt.
$$
If $t_0 \in I$, we also split $S(u) = S_{\leq t_0}(u) + S_{\geq t_0}(u)$, where
$$
S_{\leq t_0}(u) := \int_{I \cap (-\infty,t_0]} \int_{\R^d} |u(t,x)|^{2(d+2)/d}\ dx dt
$$
and
$$
S_{\geq t_0}(u) := \int_{I \cap [t_0,+\infty)} \int_{\R^d} |u(t,x)|^{2(d+2)/d}\ dx dt.
$$

We shall use this scattering size to measure blowup and scattering of solutions:

\begin{definition}[Blowup]
We say that a solution $u: I \times \R^d \to \C$ \emph{blows up forward in time} if $S_{\geq t}(u) = +\infty$ for all $t$ in the
interior of $I$ (or equivalently, for at least one such $t$).  Similarly, we say that $u$ \emph{blows up backward in time} if $S_{\leq t}(u) = +\infty$
for all (or one) $t$ in the interior of $I$.  Note that blowup may occur either at finite or infinite time, depending on whether the relevant
endpoint of $I$ is finite or infinite.
\end{definition}

\begin{definition}[Scattering]
Let $u: I \times \R^d \to \C$ be a solution.  If $u_+ \in L^2_x(\R^d)$, we say that $u$ \emph{scatters forward in time} to
$e^{it\Delta} u_+$ if $\sup I = +\infty$ and $\lim_{t \to +\infty} M( u(t) - e^{it\Delta} u_+ ) = 0$. In particular, this implies $M(u) = M(u_+)$.
Similarly, if $u_- \in L^2_x(\R^d)$, we say that $u$ \emph{scatters backward in time} to $e^{it\Delta} u_-$ if $\inf I = -\infty$ and
$\lim_{t \to -\infty} M( u(t) - e^{it\Delta} u_- ) = 0$.
\end{definition}

The use of the scattering size to measure blowup is justified by the following standard well-posedness theorem for this equation:

\begin{theorem}[Well-posedness]\label{lwp}  Let $u_0 \in L^2_x(\R^d)$ and $t_0 \in \R$.  Then there exists a unique maximal-lifespan solution
$u: I \times \R^d \to \C$ to \eqref{nls} with $t_0 \in I$ and initial data $u(t_0)=u_0$.  Furthermore:
\begin{itemize}
\item[(i)] (Local existence) The interval $I$ is an open subset of $\R$.
\item[(ii)] (Mass conservation) We have $M(u(t)) = M(u_0)$ for all $t \in I$; thus, we may define $M(u) := M(u(t))$.
\item[(iii)] (Forward scattering) If $u$ does not blow up forward in time, then $\sup I = +\infty$, and furthermore $u$ scatters forward in time to
$e^{it\Delta} u_+$ for some $u_+ \in L^2_x(\R^d)$.  Conversely, if $u_+ \in L^2_x(\R^d)$, then there exists a unique maximal-lifespan solution $u$
which scatters forward in time to $e^{it\Delta} u_+$.
\item[(iv)] (Backward scattering) If $u$ does not blow up backward in time, then $\inf I = -\infty$, and furthermore $u$ scatters backward in time to
$e^{it\Delta} u_-$ for some $u_- \in L^2_x(\R^d)$.  Conversely, if $u_- \in L^2_x(\R^d)$, then there exists a unique maximal-lifespan solution $u$
which scatters backward in time to $e^{it\Delta} u_-$.
\item[(v)] (Small data scattering) There exists a constant $C_d > 0$ depending only on dimension such that if $M(u) \leq C_d^{-1}$,
then $S(u) \leq C_d M(u)^{(d+2)/d}$.  In particular, no blowup occurs and we have global existence and scattering in both directions.
\item[(vi)] (Uniformly continuous dependence) For every $A > 0$ and $\eps > 0$ there exists $\delta > 0$ with the following property:
if $u: I \times \R^d \to \C$ is a solution (not necessarily maximal-lifespan) such that $M(u), S(u) \leq A$ and $t_0 \in I$, $v_0 \in L^2_x(\R^d)$
are such that $M(u(t_0)-v_0) \leq \delta$, then there exists a solution $v: I \times \R^d \to \C$ with $v(t_0) = v_0$ such that $S(u-v) \leq \eps$
and $M(u(t)-v(t)) \leq \eps$ for all $t \in I$.
\end{itemize}
\end{theorem}

\begin{proof} See  \cite{caz}, \cite{cwI}, or \cite{tsutsumi}.  The arguments rely primarily on the Strichartz estimate \eqref{strich} below and on the
H\"older inequality
\begin{equation}\label{holder}
 \begin{split}
 \| F(u) - F(v) &\|_{L^{2(d+2)/(d+4)}_{t,x}(I \times \R^d)} \leq C'_d \| u-v\|_{L^{2(d+2)/d}_{t,x}(I \times \R^d)}\\
\times &\left[ \|u \|_{L^{2(d+2)/d}_{t,x}(I \times \R^d)} + \|v\|_{L^{2(d+2)/d}_{t,x}(I \times \R^d)}\right]^{4/d},
\end{split}
\end{equation}
for some constant $C'_d$ depending only on $d$.  Thus, for instance, the constant $C_d$ in (v) is a simple algebraic combination of $C'_d$ and
the constant $C''_d$ appearing in \eqref{strich}.  To establish (vi) when $S(u)$ is large, one first subdivides the time interval $I$ into a
finite number of subintervals $J$, on each of which the scattering size is small, and then applies the Strichartz estimate \eqref{strich} on each
subinterval $J$ in turn.  While the scattering results and the local existence results can be proven separately via the Strichartz estimates,
they can also be deduced from each other using the pseudoconformal symmetry
\begin{equation}\label{pc}
u(t,x) \mapsto \frac{1}{|t|^{d/2}} e^{i|x|^2/4t} u(\frac{-1}{t}, \frac{x}{t})
\end{equation}
(see e.g. \cite{blue}).
\end{proof}

We now investigate the precise relationship between the mass $M(u)$ and the scattering size $S(u)$.  For any mass $m \geq 0$, let $A(m)$ denote the
quantity
\begin{equation}\label{Adef}
 A(m) := \sup \{ S(u): M(u) \leq m \},
\end{equation}
where the supremum is taken over all solutions of mass at most $m$. Thus, $A: [0,+\infty) \to [0,+\infty]$ is a monotone non-decreasing function of
$m$. From Theorem~\ref{lwp}(v) we have
\begin{equation}\label{amd}
 A(m) \leq C_d m^{(d+2)/d} \hbox{ for } m \leq C_d^{-1}.
\end{equation}
In particular, $A$ is finite for small $m$.  On the other hand, from Theorem \ref{lwp}(vi) we see that $A$ is left-continuous. Thus, there must exist
a unique \emph{critical mass} $0 < m_0 = m_0(\mu,d) \leq +\infty$ such that $A(m)$ is finite for all $m < m_0$ but infinite for all $m \geq m_0$.  By
construction, we see that
$$ S(u) \leq A(M(u))$$
for any solution $u$. In particular, if $u$ scatters forward or backward in time to $e^{it\Delta} u_\pm$, then
$$ S(u) \leq A(M(u_\pm)).$$
From Theorem \ref{lwp} we thus see that one has global well-posedness and scattering in $L^2_x(\R^d)$ whenever the mass is strictly less than the critical mass $m_0$.

It is thus of interest\footnote{The more quantitative question concerning the determination of the order of magnitude of $A(\cdot)$ is also an
extremely interesting problem, which we will not address here; any argument which goes through the compactness method will give either
a very terrible bound on $A(\cdot)$, or no bound at all.} to determine the critical mass $m_0$.  In the focusing case $\mu = -1$, it is known that
$m_0$ is finite.  Indeed, if $Q: \R^d \to \R^+$ is the ground state, that is, the unique radial positive Schwartz solution to the elliptic equation
\begin{equation}\label{ground}
 \Delta Q + Q^{1+4/d} = Q,
\end{equation}
then we have the explicit maximal-lifespan solution $u(t,x) = e^{it} Q(x)$ for $t \in \R$ which blows up both forward and backwards in
time\footnote{The forward and backward blowup times here are both infinite, but one can move the forward blowup time (say) to be finite by using the
pseudoconformal transform, though by doing so one removes the backward blowup.  See \cite{keraani}, \cite{tao-lens} for further discussion of these
different types of blowup and the relationship with the pseudoconformal transform.}, and so we have $m_0 \leq M(Q)$.

On the other hand, in the defocusing case $\mu=+1$ there is no analogous ground state.  This leads to

\begin{conjecture}[Scattering conjecture]\label{conj}
In the defocusing case ($\mu=+1$) we have $m_0=+\infty$, while in the focusing case ($\mu=-1$) we have $m_0 = M(Q)$.
\end{conjecture}

\begin{remark}
This conjecture is related to a number of other statements concerning nonlinear Schr\"odinger and generalized Korteweg-de Vries equations; see
\cite{tao-lens}, \cite{tao-gkdv}, \cite{tvz}.  For solutions in the energy class $H^1_x(\R^d)$ it is known that one has global existence for all
masses in the defocusing case and for masses $M(u) < M(Q)$ in the focusing case; see \cite{merle}, \cite{weinstein}, although control of the
scattering size $S(u)$ has not yet been obtained in these cases (in particular, blowup at infinite time has not yet been ruled out for these
solutions). There has been some progress in lowering the regularity of $H^1_x(\R^d)$ for these results, see \cite{bourg.2d}, \cite{ckstt:cubic},
\cite{cr}, \cite{hmidi}, \cite{hmidi2}, \cite{keraani}, \cite{vz}, but these methods are unlikely to reach all the way to the scale-invariant
regularity $L^2_x(\R^d)$.  At this regularity the global well-posedness and scattering problems are in fact equivalent (see \cite{begout}, \cite{blue},
\cite{keraani}, \cite{tao-lens}).
\end{remark}

We will not prove Conjecture \ref{conj} here.  However, we shall establish a basic first step towards this conjecture, which is to reduce matters to
understanding a very special subclass of solutions, namely those solutions which are \emph{almost periodic modulo the symmetries} of phase rotation,
modulation, spatial translation, and scaling.  Note that the last three symmetries are non-compact; this triple failure of compactness is a major
source of difficulty in analysing this equation.  For comparison, the energy-critical NLS (with $F(z) = \mu |z|^{4/(d-2)} z$) only has non-compactness
arising from spatial translation and scaling symmetry, while for subcritical NLS there is only the non-compactness from spatial translation.
If one assumes spherical symmetry, then one can eliminate the modulation and spatial translation sources of non-compactness, leaving only scaling
if the equation is critical.

We now pause to describe these symmetries more formally.

\subsection{The symmetry group $G$}

\begin{definition}[Symmetry group]
For any phase $\theta \in \R/2\pi \Z$, position $x_0 \in \R^d$, frequency $\xi_0 \in \R^d$, and scaling parameter $\lambda > 0$, we define the
unitary transformation $g_{\theta,x_0,\xi_0,\lambda}: L^2_x(\R^d) \to L^2_x(\R^d)$ by the formula
$$
g_{\theta, \xi_0, x_0, \lambda} f(x) :=  \frac{1}{\lambda^{d/2}} e^{i\theta} e^{i x \cdot \xi_0 } f( \frac{x-x_0}{\lambda} ).
$$
We let $G$ be the collection of such transformations\footnote{There are other symmetries one could add here, such as the rotations and the quadratic
modulations $f(x) \mapsto e^{i \tau |x|^2} f(x)$, in order to incorporate the pseudoconformal symmetry \eqref{pc}, but we will not need them here. We
also avoid the conjugation symmetry $f(x) \mapsto \overline{f(x)}$ as this reverses the arrow of time.}; this is a group with identity $g_{0,0,0,1}$,
inverse $g_{\theta,\xi_0,x_0,\lambda}^{-1} = g_{-\theta - x_0 \xi_0,\,-\lambda\xi_0,\,-x_0/\lambda,\,\lambda^{-1}}$ and group law
$$
g_{\theta, \xi_0, x_0, \lambda} g_{\theta', \xi'_0, x'_0, \lambda'} = g_{ \theta + \theta' - x_0 \cdot \xi'_0/\lambda ,\, \xi_0 + \xi'_0 / \lambda,\,
x_0 + \lambda x'_0 ,\, \lambda \lambda' }.
$$
Note that we have the factorisation
$$g_{\theta, \xi_0, x_0, \lambda} = g_{\theta, 0,0,1} g_{0,\xi_0,0,1} g_{0,0,x_0,1} g_{0,0,0,\lambda},$$
and so $G$ is generated by phase rotations, frequency modulations, translations, and dilations. We let $G\backslash L^2_x(\R^d)$ be the moduli space
of $G$-orbits $Gf := \{gf: g \in G\}$ of $L^2_x(\R^d)$, endowed with the usual quotient topology. If $u: I \times \R^d \to \C$ is a function, we
define $T_{g_{\theta,\xi_0,x_0,\lambda}} u: \lambda^2 I \times \R^d \to \C$ where $\lambda^2 I := \{ \lambda^2 t: t \in I \}$ by the formula
$$
(T_{g_{\theta, \xi_0, x_0, \lambda}} u)(t,x) :=  \frac{1}{\lambda^{d/2}} e^{i\theta} e^{i x \cdot \xi_0 } e^{-it|\xi_0|^2}
u( \frac{t}{\lambda^2}, \frac{x-x_0 - 2\xi_0 t}{\lambda} )$$
or equivalently
$$
(T_{g_{\theta, \xi_0, x_0, \lambda}} u)(t) = g_{\theta - t |\xi_0|^2, \xi_0, x_0 + 2 \xi_0 t, \lambda}(u(\frac{t}{\lambda^2}))
$$
Observe that the map $g \mapsto T_g$ is a group action of $G$.
\end{definition}

\begin{remark}[Invariances]
As $G$ is a unitary group we have $M( g u_0 ) = M(u_0)$ for any $u_0 \in L^2_x(\R^d)$.  The NLS \eqref{nls} is invariant under phase rotations,
Galilean transforms, and spatial translations, and so we see that the group action $g \mapsto T_g$ maps solutions to solutions
(and also preserves solutions $e^{it\Delta} u_0$ to the linear equation).  Furthermore, if $u$ is a solution and $g \in G$, we see that $M(T_g u) =
M(u)$ and $S(T_gu) = S(u)$, and thus the action preserves both the mass and the scattering size. Because of these symmetries, the evolution of NLS
not only foliates $L^2_x(\R^d)$ into curves $\{ u(t): t \in I\}$, where $u: I \times \R^d \to \C$ is a maximal-lifespan solution, but also foliates
the moduli space $G\backslash L^2_x(\R^d)$ into curves $\{ G u(t): t \in I\}$, though the latter curves are only parameterized affinely\footnote{In
other words, there is no single canonical time parameterisation of these curves; instead, there is an equivalence class under affine transformations
$t \mapsto at+b$ of such parameterisations.  If one were to enlarge $G$ by adding quadratic modulations $f \mapsto e^{i\tau |x|^2} f$, thus allowing
the introduction of the pseudoconformal symmetry \eqref{pc}, then the resulting curves in the moduli space would only be parameterized projectively.
Actually, in this case it would be sensible to use a lens-transformed time variable $\tan^{-1} t$, as in \cite{tao-lens}, in order to eliminate
artificial coordinate singularities.}, due to the fact that the action of $G$ on solutions rescales the time variable.  As it turns out, this
quotienting out by $G$ will serve to compactify the dynamics of NLS for certain ``minimal-mass blowup solutions'' which we shall consider shortly.
\end{remark}

\begin{remark}[Topology of $G$]
If we give the unitary operators in $G$ the strong or weak operator topology then the identification of $G$ with $\R/2\pi\Z \times \R^d \times \R^d
\times (0,+\infty)$ is easily seen to be a homeomorphism, and $G$ now has the structure of a $2d+2$-dimensional Lie group. In the strong operator
topology, $G$ is closed in the space $B(L^2_x(\R^d))$ of bounded linear operators on $L^2_x(\R^d)$, but in the weak operator topology, $G$ has $0$ as
an adherent point, and $G \cup \{0\}$ is homeomorphic to the one-point compactification of the cylinder $\R/2\pi\Z \times \R^d \times \R^d \times
(0,+\infty)$.  In particular, if $g_n \in G$ does not converge to zero in the weak operator topology, then it has a subsequence which converges in
the strong operator norm topology.  In other words, $G$ is a group of \emph{dislocations} in the sense of \cite{tintarev}, and is thus a suitable
group for constructing a concentration-compactness theory.
\end{remark}

\begin{remark}[Topology of $G \backslash L^2_x(\R^d)$]\label{topg}
From the fact that $G$ is a group of dislocations one can verify that the orbits of $G$ in $L^2_x(\R^d)$ are closed in the strong $L^2_x(\R^d)$
topology (the zero orbit has to be treated separately).  From this (and the linearity of the group action) we see that any two orbits of $G$ must be
a non-zero distance apart.  Thus $G\backslash L^2_x(\R^d)$ is in fact a metric space.  Indeed, one can show (again using the dislocation property)
that any convergent (or Cauchy) sequence in $G \backslash L^2_x(\R^d)$ is the projection of a convergent (or Cauchy) sequence in $L^2_x(\R^d)$.  In
particular, since $L^2_x(\R^d)$ is a complete metric space, $G\backslash L^2_x(\R^d)$ is also.
\end{remark}

\subsection{Main result}

Let us say that a function $u \in C^0_{t,\loc} L^2_x(I \times \R^d)$ is \emph{almost periodic modulo $G$} if the quotiented orbit $\{ Gu(t): t \in I
\}$ is a precompact subset of $G\backslash L^2_x(\R^d)$ (i.e. its closure is compact, thus every sequence in the orbit has a convergent subsequence
in $G\backslash L^2_x(\R^d)$).  Equivalently (by Remark \ref{topg}), $u$ is almost periodic modulo $G$ if there exists a compact subset $K$ of
$L^2_x(\R^d)$ such that $u(t) \in GK$ for all $t \in I$, or in other words, $g(t)^{-1} u(t) \in K$ for all $t \in I$ and some function $g: I \to
G$.  For instance, if one considered a (hypothetical) periodic ``breather'' solution to NLS and applied a Galilean transform to it, the resulting
``travelling breather'' solution would not be periodic or almost periodic in $L^2_x(\R^d)$ in the classical sense, but would be periodic (and thus
almost periodic) modulo the symmetry group $G$ (in fact, modulo $G$ the travelling breather traverses the same orbit as the stationary breather).

We can now state the main result of this paper.

\begin{theorem}[Reduction to almost periodic solutions]\label{main}
Fix $\mu$ and $d$, and suppose that the critical mass $m_0$ is finite.  Then there exists a maximal-lifespan solution $u$ of mass
exactly $m_0$ which blows up both forward and backward in time.  Furthermore, any maximal-lifespan solution of mass $m_0$ which blows up
both forward and backward in time is almost periodic modulo $G$.
\end{theorem}

\begin{example}\label{blowex}
Suppose that Conjecture \ref{conj} is true in the focusing case, so $m_0 = M(Q)$.  Then the soliton solution $e^{it} Q(x)$ (or more generally, the
other soliton solutions which range in the orbit $GQ$) will be a solution of the type claimed in Theorem \ref{main}.  On the other hand, if one
applies a pseudoconformal transformation \eqref{pc} to a soliton solution, then one obtains two solutions which blow up in only one time direction,
and one will only have almost periodicity in that direction and not in the direction without blowup\footnote{One can recover almost periodicity in
both directions by enlarging the symmetry group $G$ to contain quadratic modulations $f \mapsto e^{i\tau |x|^2} f$, but this complicates the role of
the time parameter and we shall avoid doing this here.}.  Indeed one can easily show that a maximal-lifespan solution which is almost periodic modulo
$G$ and not identically zero \emph{must} blow up both forward and backward in time; for if it did not blow up forward in time (say), then it scatters
to a linear solution $e^{it\Delta} u_+$ by Theorem \ref{lwp}(iii), which can easily be seen to be incompatible with almost periodicity modulo $G$ by
obtaining stationary phase asymptotics for the linear solution (or by using the pseudoconformal transform \eqref{pc}).
\end{example}

\begin{remark}
This result is essentially implicit in the work of Keraani, \cite{keraani}, in dimensions $d=1,2$, and in principle follows in the general dimension
case by using the results of Begout-Vargas, \cite{begout}, (which we also use in this paper), although this is not stated explicitly.  See also the
analysis of Merle and Vega, \cite{mv}, concerning blowup solutions (with mass possibly larger than $m_0$).  Indeed, our methods are a
combination of the concentration-compactness arguments in \cite{begout}, \cite{mv}, \cite{keraani} and the induction on energy
method in \cite{bourg.critical}, \cite{gopher}, \cite{rv}, \cite{visan}.  The high degree of compatibility between these two
arguments was first observed in \cite{merlekenig}, where an analogue of Theorem~\ref{main} was proven for the energy-critical NLS; our
approach is in fact very similar to that in \cite{merlekenig}.  Indeed, it seems that the phenomenon of existence of blowup
solutions at critical levels of a conserved quantity, which are almost periodic modulo the symmetries of the problem, is a very general one,
essentially being a consequence of a sufficiently strong stability and concentration-compactness theory for such equations.
\end{remark}

\begin{remark}
In view of Theorem \ref{main}, we see that to prove Conjecture~\ref{conj} it suffices to show that in the defocusing case there do not exist any
maximal-lifespan solutions $u$ which are almost periodic modulo $G$ other than the zero solution, while in the focusing case it suffices to establish
the same statement under the additional hypothesis $M(u) < M(Q)$.  In the language of Martel and Merle, \cite{martel}, this reduces matters to
establishing a ``Liouville theorem'' for these equations.
\end{remark}

One can phrase the property of almost periodicity modulo $G$ in a more ``quantitative'' sense (in the spirit of \cite{gopher}, \cite{rv},
\cite{visan}) as follows:

\begin{lemma}
Let $u \in C^0_{t,\loc} L^2_x(I \times \R^d)$ be almost periodic modulo $G$.  Then there exist functions $x: I \to \R^d$, $\xi: I \to \R^d$, and $N:
I \to (0,+\infty)$ with the property that for every $\eta > 0$ there exists $0 < C(\eta) < \infty$ such that we have the spatial concentration
estimate
\begin{equation}\label{spaceconc}
\int_{|x - x(t)| \geq C(\eta)/N(t)} |u(t,x)|^2\ dx \leq \eta
\end{equation}
and frequency concentration estimate
\begin{equation}\label{freqconc}
\int_{|\xi - \xi(t)| \geq C(\eta) N(t)} |\hat u(t,\xi)|^2\ d\xi \leq \eta
\end{equation}
for all $t \in I$.
\end{lemma}

\begin{remark}
As stated, the quantity $C(\eta)$ depends on $u$.  However, if one restricts attention to minimal-mass almost periodic solutions, that is,
$M(u) = m_0$, then one can make the quantity $C(\eta)$ independent of $u$, by arguing by contradiction and repeating the arguments in
Section~\ref{keysec} below; we omit the details.
\end{remark}

\begin{remark}
Informally, this lemma asserts that the mass $u(t)$ is spatially concentrated in the ball $\{ x: x = x(t) + O( 1/N(t) ) \}$ and is frequency
concentrated in the ball $\{ \xi: \xi = \xi(t) + O(N(t)) \}$.  Note that we have currently no control as to how $x(t)$, $N(t)$, $\xi(t)$ vary in
time; obtaining such control is of course very important in understanding the dynamics of these solutions, but this requires additional techniques
(e.g. monotonicity formulae, conservation laws, use of Duhamel formula) which we do not pursue here (but see \cite{tvz}).  Indeed, one can view the
results here as essentially the limit of what one can say about the minimal-mass blowup solutions to \eqref{nls} using only the perturbative theory
and the mass conservation law.
\end{remark}

\begin{proof}
By hypothesis, $u(t)$ lies in $GK$ for some compact subset $K$ in $L^2_x(\R^d)$.  A simple compactness argument shows that for every $\eta > 0$ there
exists $0 < C(\eta) < \infty$ (depending on $K$) such that
$$ \int_{|x| \geq C(\eta)} |f(x)|^2\ dx \leq \eta$$
and frequency concentration estimate
$$ \int_{|\xi| \geq C(\eta)} |\hat f(\xi)|^2\ d\xi \leq \eta$$
for all $f \in K$.  The claim then follows by inspecting what the symmetry group $G$ does to the spatial and frequency distribution of the mass of a
function.
\end{proof}

Combining this with Theorem \ref{main} we obtain

\begin{corollary}\label{cormain} Fix $\mu$ and $d$, and suppose that $m_0$ is finite.
Then there exists a maximal-lifespan solution $u \in C^0_{t,\loc} L^2_x(I \times \R^d)$ of mass exactly $m_0$ which blows up both forward and
backward in time, and functions $x: I \to \R^d$, $\xi: I \to \R^d$, and $N: I \to (0,+\infty)$ with the property that for every $\eta > 0$ there
exists $0 < C(\eta) < \infty$ (depending on $\mu,d,m_0$) such that we have the concentration estimates \eqref{spaceconc}, \eqref{freqconc} for all $t
\in I$.
\end{corollary}

\begin{remark} A very similar result appears in \cite{gopher}, \cite{rv}, \cite{visan} for the energy-critical NLS,
though there one considers solutions which ``almost blow up'' in the sense that $S(u)$ is huge rather than infinite.  As a consequence, the bounds
\eqref{spaceconc}, \eqref{freqconc} do not hold for all $\eta$, but only for all $\eta$ larger than an extremely small positive quantity.  It is
possible to repeat the induction on energy arguments in \cite{gopher}, \cite{rv}, \cite{visan}, changing the numerology appropriately, and using some
different estimates (notably the bilinear restriction estimate from \cite{tao:bilinear}) to give an alternate proof of Corollary \ref{cormain}; we
sketch this alternate derivation in Section~8. In principle, this more quantitative approach, not relying explicitly on compactness, gives some
explicit bounds on the quantity $C(\eta)$, although these bounds are extremely poor (see \cite{gopher} for some further discussion).  In any
event, the induction on energy and concentration-compactness arguments, despite many superficial differences, are in fact closely related, sharing
many of the same underlying estimates and ideas.  One difference is that in the quantitative approach in \cite{gopher}, \cite{rv},
\cite{visan}, one does not use the full power of concentration-compactness, but merely settles for extracting a single ``bubble'' of concentration.
The price one pays for this simplification is that one must then work significantly harder to show that the evolution of such bubbles are
sufficiently decoupled from the rest of the solution, for instance one needs to use tools such as approximate finite speed of propagation and
persistence of positive and negative regularities, as well as a greater reliance on the bilinear restriction estimates from \cite{tao:bilinear}
(see Section~8).
\end{remark}

In a sequel to this paper, \cite{tvz}, we shall use Corollary \ref{cormain} (or more precisely, the analogue of this corollary for spherically
symmetric solutions, see Section \ref{spherical-sec}) to establish global well-posedness and scattering for the defocusing NLS for spherically
symmetric $L^2_x(\R^d)$ solutions in dimensions $d \geq 3$.  We expect that Theorem \ref{main} will similarly be useful for the lower dimensional
case, for focusing nonlinearities, and for non-radial data, and hope to address some of these issues in future work.

\textbf{Acknowledgements}: This research was partially conducted during the period Monica Visan was employed by the Clay
Mathematics Institute as a Liftoff Fellow.  This material is based upon work supported by the National Science Foundation under agreement
No. DMS--0111298.  Any opinions, findings, and conclusions or recommendations expressed in this material are those of the authors and do not
reflect the views of the National Science Foundation.  The third author was supported by the NSF grant No. 10601060 (China).
The authors thank Sahbi Keraani, Carlos Kenig, and Frank Merle for helpful comments.

\section{The key convergence result}\label{keysec}

The proof of Theorem \ref{main} rests on the following key proposition, asserting a certain compactness modulo $G$ in blowup sequences of solutions
with mass less than or equal to the critical mass.

\begin{proposition}[Palais-Smale condition modulo $G$]\label{mainprop}  Fix $\mu$ and $d$, and suppose that $m_0$ is finite.
Let $u_n: I_n \times \R^d \to \C$ for $n=1,2,\ldots$ be a sequence of solutions and $t_n \in I_n$ a sequence of times such that $\limsup_{n \to
\infty} M(u_n)= m_0$ and
\begin{equation}\label{under}
\lim_{n \to \infty} S_{\geq t_n}(u_n) =
\lim_{n \to \infty} S_{\leq t_n}(u_n) = +\infty.
\end{equation}
Then the sequence $G u_n(t_n)$ has a subsequence which converges in the $G \backslash L^2_x(\R^d)$ topology.
\end{proposition}

\begin{remark}
The hypothesis \eqref{under} asserts that the sequence $u_n$ asymptotically blows up both forward and backward in time.  Both components of this
hypothesis are essential, as can be seen by considering the examples in Example \ref{blowex} which only blow up in one direction (and whose orbit is
non-compact in the other direction, even after quotienting out by $G$).
\end{remark}

We prove this result in Section \ref{propsec}.  For now, let us assume it and conclude the proof of Theorem \ref{main}.

\begin{proof}[Proof of Theorem \ref{main} assuming Proposition \ref{mainprop}]
By definition of $m_0$ we can find a sequence $u_n: I_n \times \R^d \to \C$ of solutions with $M(u_n) \leq m_0$ and $\lim_{n \to \infty} S(u_n) =
+\infty$. Without loss of generality we can take the $u_n$ to have maximal lifespan.  By choosing $t_n \in I_n$ to be the median time of the
$L^{2(d+2)/d}_{t,x}$ norm of $u_n$ (cf. the ``middle thirds'' trick in \cite{bourg.critical}, \cite{gopher}) we can thus arrange that \eqref{under}
holds.  By time translation invariance we may take $t_n = 0$.  We then apply Proposition \ref{mainprop}, and after passing to a subsequence if
necessary, we can locate $u_0 \in L^2_x(\R^d)$ such that $G u_n(0)$ converges in the $G\backslash L^2_x(\R^d)$ topology to $Gu_0$; thus, there are
group elements $g_n \in G$ such that $g_n u_n(0)$ converges strongly in $L^2_x(\R^d)$ to $u_0$.  By applying the group action $T_{g_n}$ to the
solutions $u_n$ we may take the $g_n$ to all be the identity, thus $u_n(0)$ now converges strongly in $L^2_x(\R^d)$ to $u_0$. In particular this
implies $M(u_0) \leq m_0$.

Let $u: I \times \R^n \to \C$ be the maximal-lifespan solution with initial data $u(0) = u_0$ as given by Theorem \ref{lwp}.  We claim that $u$ blows
up both forward and backward in time.  Indeed, if $u$ does not blow up forward in time (say), then by Theorem \ref{lwp}(iii) we have $[0,+\infty)
\subset I$ and $S_{\geq 0}(u) < \infty$.  By Theorem \ref{lwp}(vi), this implies for sufficiently large $n$ that $[0,+\infty) \subset I_n$ and
$$ \limsup_{n \to \infty} S_{\geq 0}( u_n ) < \infty,$$
contradicting \eqref{under}.  Similarly if $u$ blows up backward in time.  By the definition of $m_0$ this forces $M(u_0) \geq m_0$, and hence
$M(u_0)$ must be exactly $m_0$.

It remains to show that a solution $u$ which blows up both forward and backward in time is almost periodic modulo $G$.  Consider an arbitrary
sequence $G u(t'_n)$ in $\{ G u(t): t \in I \}$. Now, since $u$ blows up both forward and backward in time, but is locally in $L^{2(d+2)/d}_{t,x}$,
we have
$$ S_{\geq t'_n}(u) = S_{\leq t'_n}(u) = \infty.$$
Applying Proposition \ref{mainprop} once again we see that $G u(t'_n)$ does have a convergent sequence in $G \backslash L^2_x(\R^d)$.  Thus, the
orbit $\{ G u(t): t \in I \}$ is precompact in $G \backslash L^2_x(\R^d)$ as desired.
\end{proof}

\begin{remark}
One can modify the above argument to show that if a solution $u: I \times \R^d \to \C$ with mass $m_0$ blows up forward in time only, then the
restriction of $u$ to any subinterval $I_{\geq t}$ for $t \in I$ will be almost periodic, thus one has almost periodicity forward in time only (cf.
\cite[Proposition 4.2]{merlekenig}).  We omit the details.
\end{remark}

It remains to prove Proposition \ref{mainprop}.  In order to do so we need to recall two standard tools, namely a stability result for NLS and a
concentration-compactness result for solutions to the linear equation.  This is the purpose of the next two sections.

\section{A stability lemma}\label{stab-sec}

We have the following standard Strichartz estimate (\cite{strich}, \cite{gv}):
if $u: I \times \R^d \to \C$ solves the inhomogeneous Schr\"odinger equation
$$ iu_t + \Delta u = F; \quad u(t_0) = u_0$$
for some $t_0 \in I$, $u_0 \in L^2_x(\R^d)$, and $F \in L^{2(d+2)/(d+4)}_{t,x}(I \times \R^d)$ in the integral (Duhamel) sense, that is,
$$ u(t) = e^{i(t-t_0)\Delta} u_0 - i \int_{t_0}^t e^{i(t-t')\Delta} F(t')\ dt',$$
then we have
\begin{equation}\label{strich}
\| u \|_{C^0_t L^2_x(I \times \R^d)} + \| u \|_{L^{2(d+2)/d}_{t,x}(I \times \R^d)} \leq C''_d ( \|u_0\|_{L^2_x(\R^d)} +
\|F\|_{L^{2(d+2)/(d+4)}_{t,x}(I \times \R^d)})
\end{equation}
for some constant $0 < C''_d < \infty$ depending only on the dimension $d$.  Many more Strichartz estimates are available (see e.g. \cite{tao:keel}),
but this is the only one we shall need here.  (A bilinear refinement of \eqref{strich} will however be implicit in the proof of Theorem \ref{lp}
below.)

As remarked previously, this estimate underlies all the results in Theorem \ref{lwp}, which we rely on extensively in this paper.  We shall also need
a variant of this theorem in which one starts with an \emph{approximate} solution to \eqref{nls} and perturbs it to an exact solution:

\begin{lemma}[Stability, \cite{tvz}]\label{stab} Fix $\mu$ and $d$.
For every $A > 0$ and $\eps > 0$ there exists $\delta > 0$ with the following property: if $u: I \times \R^d \to \C$ is such that $S(u) \leq A$ and
that $u$ approximately solves \eqref{nls} in the sense that
\begin{equation}\label{ufd}
 \| iu_t + \Delta u - F(u) \|_{L^{2(d+2)/(d+4)}_{t,x}(I \times \R^d)} \leq \delta,
\end{equation}
and $t_0 \in I$, $v_0 \in L^2_x(\R^d)$ are such that $S_I(e^{i(t-t_0)\Delta} (u(t_0)-v_0)) \leq \delta^{2(d+2)/d}$, then there exists a solution $v:
I \times \R^d \to \C$ to \eqref{nls}with $v(t_0) = v_0$ such that $S(u-v) \leq \eps$.
\end{lemma}

\begin{remark}
This generalizes Theorem \ref{lwp}(vi) (except that we no longer control the mass of $u(t)-v(t)$), because we now allow $iu_t + \Delta u - F(u)$ to
be small but nonzero.  It also implies the existence and uniqueness of maximal-lifespan solutions in Theorem \ref{lwp}.  Interestingly, the masses of
$u$ and $v_0$ do not directly appear in this lemma, though it is necessary that these masses are finite. Analogous stability results for the
energy-critical NLS (in $\dot H^1(\R^d)$ instead of $L^2_x(\R^d)$, of course) have appeared in \cite{gopher}, \cite{merlekenig}, \cite{rv},
\cite{TaoVisan}, \cite{visan}.  The mass-critical case is in fact slightly simpler as one does not need to deal with the presence of a
derivative in the regularity class.
\end{remark}

\begin{proof}  (Sketch)
Let us first establish the claim when $A$ is sufficiently small depending on $d$. Let $v: I' \times \R^d \to \C$ be the maximal-lifespan solution
with initial data $v(t_0) = v_0$.  Writing $v = u+w$ on the interval $I'' := I \cap I'$, we see that
$$ iw_t + \Delta w = (F(u+w) - F(u)) - (iu_t + \Delta u - F(u))$$
and $S_{I''}(e^{i(t-t_0)\Delta}w(t_0) \leq \delta^{2(d+2)/d}$.  Thus, if we set
$$X := \| w \|_{L^{2(d+2)/d}_{t,x}(I'' \times \R^d)}, $$
by \eqref{strich}, \eqref{ufd}, and the triangle inequality,  we have
$$ X \leq \delta + C''_d [ \| F(u+w) - F(u) \|_{L^{2(d+2)/(d+4)}_{t,x}(I'' \times \R^d)} + \delta ],$$
and hence, by \eqref{holder} and the hypothesis $S(u) \leq A$,
$$ X \leq \tilde C_d [A^{2/(d+2)} X + X^{1+4/d} + \delta ]$$
where $\tilde C_d$ depends only on $d$. If $A$ is sufficiently small depending on $d$, and $\delta$ is sufficiently small depending on $\eps$ and
$d$, then standard continuity arguments give $X \leq \eps$ as desired.  To handle the case when $A$ is large, simply iterate the case when $A$ is
small (shrinking $\delta$, $\eps$ repeatedly) after a subdivision of the time interval $I$.
\end{proof}

\begin{remark} From \eqref{strich}, we see that we can replace the hypothesis
 $S( e^{i(t-t_0)\Delta} ( u(t_0)-v_0) ) \leq \delta^{2(d+2)/d}$ by the essentially stronger hypothesis $M( u(t_0)-v_0 ) \leq \delta^2$, if desired.
\end{remark}

\section{Concentration-compactness}

We now recall a key concentration-compactness result of Begout-Vargas, \cite{begout}, regarding the defect of compactness in \eqref{strich}, based in
turn on earlier work of Merle-Vega, \cite{mv}, and Carles-Keraani, \cite{ck}, who handled the cases $d=2$ and $d=1$ respectively.  Because of the
time translation invariance of the linear Schr\"odinger equation, we will need to enlarge the group $G$ to contain the linear propagators
$e^{it_0\Delta}$ (though we will later be able to eventually descend back to the original group).

\begin{definition}[Enlarged group]
For any phase $\theta \in \R/2\pi \Z$, position $x_0 \in \R^d$, frequency $\xi_0 \in \R^d$, scaling parameter $\lambda > 0$, and time $t_0$, we
define the unitary transformation $g_{\theta,x_0,\xi_0,\lambda,t_0}: L^2_x(\R^d) \to L^2_x(\R^d)$ by the formula
$$g_{\theta, \xi_0, x_0, \lambda,t_0} = g_{\theta, \xi_0, x_0, \lambda} e^{it_0\Delta},$$
or in other words\footnote{Our notation here differs slightly from that in \cite{keraani}, \cite{begout}, as it is convenient for us to place the
free propagator $e^{it_0\Delta}$ at the right, rather than in the middle.  This explains for instance the discrepancy between \eqref{div} and
\cite[Definition 5.3]{begout}.}
$$ g_{\theta, \xi_0, x_0, \lambda,t_0} f(x) :=  \frac{1}{\lambda^{d/2}} e^{i\theta} e^{i x \cdot \xi_0 } [e^{it_0\Delta} f]( \frac{x-x_0}{\lambda} ).$$
Let $G'$ be the collection of such transformations.  We also let $G'$ act on global spacetime functions $u: \R \times \R^d \to \C$ by defining
$$ T_{g_{\theta, \xi_0, x_0, \lambda,t_0}} u(t,x) :=
\frac{1}{\lambda^{d/2}} e^{i\theta} e^{i x \cdot \xi_0 } e^{-it|\xi_0|^2} (e^{it_0\Delta} u)( \frac{t}{\lambda^2}, \frac{x-x_0 - 2\xi_0 t}{\lambda}),
$$
or equivalently
$$
(T_{g_{\theta, \xi_0, x_0, \lambda,t_0}} u)(t) = g_{\theta - t |\xi_0|^2, \xi_0, x_0 + 2 \xi_0 t, \lambda,t_0}(u(\frac{t}{\lambda^2})).
$$
\end{definition}

One can verify that $G'$ is indeed a group (this is basically due to the phase rotation, Galilean, translation, and scaling symmetries of the
\emph{linear} Schr\"odinger equation).  The action of $G'$ on global spacetime functions preserves solutions of the linear equation as well as the
scattering size $S(\cdot)$, but not solutions of the nonlinear equation (to achieve the latter, one would have to replace the linear propagators $e^{it_0
\Delta}$ with their nonlinear counterparts). Using either the strong or weak operator topology as before, we can identify $G'$ topologically with
$\R/2\pi \Z \times \R^d \times \R^d \times (0,+\infty) \times \R$, giving $G'$ the structure of a $2d+3$-dimensional Lie group.  Given any two
sequences $g_n, g'_n$ in $G'$, we say that $g_n$ and $g'_n$ are \emph{asymptotically orthogonal} if $(g_n)^{-1} g'_n$ diverges to infinity in $G'$
(i.e. it leaves any compact subset of $G'$ for sufficiently large $n$).  If we write explicitly
$$ g_n = g_{\theta_n, \xi_n, x_n, \lambda_n, t_n}; \quad g'_n = g_{\theta'_n, \xi'_n, x'_n, \lambda'_n, t'_n},$$
then this asymptotic orthogonality is equivalent to
\begin{equation}\label{div}
 \lim_{n \to \infty} \frac{\lambda_n}{\lambda'_n} + \frac{\lambda'_n}{\lambda_n} + |t_n \lambda^2_n - t'_n (\lambda'_n)^2|
+ |\xi_n - \xi'_n| + |x_n - x'_n| = +\infty.
\end{equation}

The terminology ``asymptotic orthogonality'' is justified by the following easy observation (a variant of the Riemann-Lebesgue lemma): if $g_n, g'_n$
are asymptotically orthogonal, then
$$ \lim_{n \to \infty} \langle g_n f, g'_n f' \rangle_{L^2_x(\R^d)} = 0 \hbox{ for all } f, f' \in L^2_x(\R^d).$$
This is essentially the assertion that $G'$ is a group of dislocations.  A variant of this is that if $v, v' \in L^{2(d+2)/d}_{t,x}(\R \times \R^d)$,
then
$$ \lim_{n \to \infty} \| |T_{g_n} v|^{1/2} |T_{g'_n} v'|^{1/2} \|_{L^{2(d+2)/d}_{t,x}(\R \times \R^d)} = 0$$
(see \cite{bg}, \cite{mv}).  From H\"older's inequality we also deduce the more general version
\begin{equation}\label{ntg}
 \lim_{n \to \infty} \| |T_{g_n} v|^{1-\theta} |T_{g'_n} v'|^\theta \|_{L^{2(d+2)/d}_{t,x}(\R \times \R^d)} = 0
 \end{equation}
for any $0 < \theta < 1$.

As a consequence of these estimates, we see that if
$g^{(j)}_n$ are sequences for $j=1,\ldots,l$ which are pairwise asymptotically orthogonal, and $f^{(1)},\ldots,f^{(l)} \in L^2_x(\R^d)$, then
$$ \lim_{n \to \infty} M( \sum_{j=1}^l g^{(j)}_n f^{(j)} ) = \sum_{j=1}^l M( f^{(j)} ).$$
Similarly, if $v^{(1)}, \ldots, v^{(l)} \in L^{2(d+2)/d}_{t,x}(\R \times \R^d)$, then we have
\begin{equation}\label{sjl}
 \lim_{n \to \infty} S( \sum_{j=1}^l g^{(j)}_n v^{(j)} ) = \sum_{j=1}^l S( v^{(j)} )
\end{equation}
(cf. \cite[Lemma 5.5]{begout}); this follows from \eqref{ntg} by using the elementary estimate
$$ \left| |\sum_{j=1}^l z_j| - |\sum_{j=1}^l |z_j|^{2(d+2)/d}|^{d/2(d+2)} \right|
\leq C_{d,l} \sum_{j \neq j'} |z_j|^{\theta} |z_{j'}|^{1-\theta}$$
for some $C_{d,l} < \infty$ and $0 < \theta < 1$, and all complex numbers $z_1,\ldots,z_l$ (this estimate can be proven by an induction on $l$).

The estimate \eqref{strich} (in the homogeneous case $F \equiv 0$) is essentially invariant under $G'$ (i.e. it is invariant under phase rotation,
Galilean transforms, spatial translation, scaling, and time translation).  Since $G'$ is (quadruply!) non-compact, this is thus a defect of
compactness for \eqref{strich}. The main result we need from \cite{begout} asserts, roughly speaking, that this is in fact the \emph{only} source of
non-compactness for \eqref{strich}.

\begin{theorem}[Linear profiles]\label{lp}\cite[Theorem 5.4]{begout}
Fix $d$. Let $u_n$, $n=1,2,\ldots$ be a bounded sequence in $L^2_x(\R^d)$.  Then (after passing to a subsequence if necessary) there exists a family
$\phi^{(j)}$, $j=1,2,\ldots$ of functions in $L^2_x(\R^d)$ and group elements $g^{(j)}_n \in G'$ for $j,n = 1,2,\ldots$ such that we have the
decomposition
\begin{equation}\label{undecomp}
u_n = \sum_{j=1}^l g^{(j)}_n \phi^{(j)} + w^{(l)}_n
\end{equation}
for all $l=1,2,\ldots$; here, $w^{(l)}_n \in L^2_x(\R^d)$ is such that its linear evolution has asymptotically vanishing scattering size:
\begin{equation}\label{sln}
\lim_{l \to \infty} \limsup_{n \to \infty} S( e^{it\Delta} w^{(l)}_n ) = 0.
\end{equation}
Moreover, $g^{(j)}_n, g^{(j')}_n$ are asymptotically orthogonal for any $j \neq j'$, and
for any $l \geq 1$ we have the mass decoupling property
\begin{equation}\label{un}
\lim_{n \to \infty} [M( u_n ) - \sum_{j=1}^l M( \phi^{(j)} ) - M( w^{(l)}_n )] = 0.
\end{equation}
\end{theorem}

\begin{remark}
The proof of the above theorem is non-trivial, combining concentration compactness arguments with a bilinear restriction estimate from
\cite{tao:bilinear} (though one can rely on simpler bilinear restriction estimates in the cases $d=1,2$).  In Section~8, we establish the main
results of this paper without recourse to this theorem, using instead the slightly simpler inverse Strichartz theorem in \cite{tao-lens} and
repeating the quantitative space and frequency concentration analysis from \cite{gopher}.
\end{remark}

\section{Proof of Proposition \ref{mainprop}}\label{propsec}

We now prove Proposition \ref{mainprop}.  By translating $u_n$ (and $I_n$) in time, we may take $t_n = 0$ for all $n$; thus,
\begin{equation}\label{under0}
\lim_{n \to \infty} S_{\geq 0}(u_n) = \lim_{n \to \infty} S_{\leq 0}(u_n) = +\infty.
\end{equation}
We then apply Theorem \ref{lp} to the bounded sequence $u_n(0)$ (passing to a subsequence if necessary) to obtain the linear profile decompositions
\begin{equation}\label{uprof}
u_n(0) =
\sum_{j=1}^l g^{(j)}_n \phi^{(j)} + w^{(l)}_n
\end{equation}
with the stated properties.  We can factorize
$$ g^{(j)}_n = h^{(j)}_n e^{it^{(j)}_n \Delta} $$
where $t^{(j)}_n \in \R$ and $h^{(j)}_n \in G$.  By refining the subsequence once for each $j$ and using a standard diagonalisation argument, we may
assume that for each $j$ that the sequence $t^{(j)}_n$, $n=1,2,\ldots$ is convergent to some time in the compactified time interval
$[-\infty,+\infty]$. If $t^{(j)}_n$ converges to some finite time $t^{(j)} \in (-\infty,+\infty)$, we may shift $\phi^{(j)}$ by the linear propagator
$e^{it^{(j)}\Delta}$, and so we may assume without loss of generality that $t^{(j)}_n$ converges either to $-\infty$, $0$, or $+\infty$. If
$t^{(j)}_n$ converges to zero, we may absorb the error $e^{it^{(j)}_n \Delta} \phi^{(j)} - \phi^{(j)}$ into the error term $w^{(l)}_n$ (which will
not significantly affect the scattering size of the linear evolution of $w^{(l)}_n$, thanks to \eqref{strich}) and so we may assume without loss of
generality in this case that $t^{(j)}_n$ is \emph{identically} zero.

From \eqref{un} we obtain the mass decoupling
\begin{equation}\label{massbound}
 \sum_{j=1}^\infty M( \phi^{(j)} ) \leq \limsup_{n \to \infty} M( u_n(0) ) \leq m_0
\end{equation}
and in particular that
$$ \sup_j M( \phi^{(j)} ) \leq m_0.$$
Suppose first that we have
\begin{equation}\label{mj}
 \sup_j M( \phi^{(j)} ) \leq m_0 - \eps
\end{equation}
for some $\eps > 0$; we will eventually show that this leads to a contradiction.  Note that $A$ is monotone increasing and finite on the interval
$[0,m_0-\eps]$.  From this and \eqref{amd}, we see that we have the bound
\begin{equation}\label{abm}
A(m) \leq B m \hbox{ for all } 0 \leq m \leq m_0 - \eps
\end{equation}
and for some finite quantity $0 < B < +\infty$ depending on $\eps$ and $d$.

We now define a nonlinear profile $v^{(j)}: \R \times \R^d \to \C$ associated to $\phi^{(j)}$ and depending on the limiting value of $t^{(j)}_n$,
as follows:

\begin{itemize}
\item If $t^{(j)}_n$ is identically zero, we define $v^{(j)}$ to be the maximal-lifespan solution with initial data $v^{(j)}(0) = \phi^{(j)}$.
\item If $t^{(j)}_n$ converges to $+\infty$, we define $v^{(j)}$ to be the maximal-lifespan solution which scatters forward in time to $e^{it\Delta} \phi^{(j)}$.
\item If $t^{(j)}_n$ converges to $-\infty$, we define $v^{(j)}$ to be the maximal-lifespan solution which scatters backward in time to $e^{it\Delta} \phi^{(j)}$.
\end{itemize}

Since $A(m_0-\eps)$ is finite, we see from \eqref{mj}, \eqref{abm}, and Theorem \ref{lwp} that $v^{(j)}$ is defined globally in time and obeys the
estimates
$$ M( v^{(j)} ) = M( \phi^{(j)} ) \leq m_0 - \eps$$
and
\begin{equation}\label{sbound}
S(v^{(j)}) \leq A( M( \phi^{(j)} ) ) \leq B M( \phi^{(j)} ).
\end{equation}
We then define the approximant $u^{(l)}_n \in C^0_t L^2_x(\R \times \R^d)$ to $u_n$ for $n,l = 1,2,\ldots$ by the formula
\begin{equation}\label{ulndef}
u^{(l)}_n(t) := \sum_{j=1}^l T_{h^{(j)}_n} [v^{(j)}(\cdot + t^{(j)}_n)](t) + e^{it\Delta} w^{(l)}_n.
\end{equation}
From \eqref{sjl}, \eqref{sln}, and the triangle inequality, followed by \eqref{massbound}, \eqref{sbound}, we have
\begin{equation}\label{sun}
\begin{split}
\lim_{l \to \infty} \lim_{n \to \infty} S( u^{(l)}_n )
&=
\lim_{l \to \infty} \sum_{j=1}^l S( v^{(j)}) \\
&\leq \lim_{l \to \infty} \sum_{j=1}^l B M(\phi^{(j)}) \\
&\leq B m_0.
\end{split}
\end{equation}

\begin{lemma}[Asymptotic agreement with initial data]  For any $l=1,2,\ldots$ we have
$$ \lim_{n \to \infty} M( u^{(l)}_n(0) - u_n(0) ) = 0.$$
\end{lemma}

\begin{proof} Comparing \eqref{ulndef} with \eqref{undecomp} and using the triangle inequality, we see that it suffices to show that
$$ \lim_{n \to \infty} M( T_{h^{(j)}_n} [v^{(j)}(\cdot + t^{(j)}_n)](0) - g^{(j)}_n \phi^{(j)} ) = 0.$$
But
\begin{align*}
M( T_{h^{(j)}_n} [v^{(j)}(\cdot + t^{(j)}_n)](0) - g^{(j)}_n \phi^{(j)} )
&= M( h^{(j)}_n [v^{(j)}(t^{(j)}_n)] - h^{(j)}_n e^{it^{(j)}_n\Delta} \phi^{(j)} ) \\
&= M( v^{(j)}(t^{(j)}_n) - e^{it^{(j)}_n\Delta} \phi^{(j)} )
\end{align*}
and the claim follows from the construction of $v^{(j)}$.
\end{proof}

\begin{lemma}[Asymptotic solvability of equation]  We have
$$ \lim_{l \to \infty} \limsup_{n \to \infty} \| (i \partial_t + \Delta) u^{(l)}_n - F(u^{(l)}_n) \|_{L^{2(d+2)/(d+4)}_{t,x}(\R \times \R^d)} = 0.$$
\end{lemma}

\begin{proof}  Write
$$ v^{(j)}_n := T_{h^{(j)}_n} [v^{(j)}(\cdot + t^{(j)}_n)].$$
By the definition of $u^{(l)}_n$, we have
$$ u^{(l)}_n = \sum_{j=1}^l v^{(j)}_n + e^{it\Delta} w^{(l)}_n$$
and
$$ (i \partial_t + \Delta) u^{(l)}_n = \sum_{j=1}^l F(v^{(j)}_n)$$
and so it suffices by the triangle inequality to show that
$$ \lim_{l \to \infty} \limsup_{n \to \infty} \| F(u^{(l)}_n - e^{it\Delta} w^{(l)}_n ) - F(u^{(l)}_n) \|_{L^{2(d+2)/(d+4)}_{t,x}(\R \times \R^d)} = 0$$
and
$$ \lim_{n \to \infty} \| F(\sum_{j=1}^l v^{(j)}_n) - \sum_{j=1}^l F(v^{(j)}_n) \|_{L^{2(d+2)/(d+4)}_{t,x}(\R \times \R^d)} = 0$$
for each $l$.

The first inequality follows immediately from \eqref{holder}, \eqref{sln}, \eqref{sun}.  For the second inequality,
we use the elementary inequality
$$ |F(\sum_{j=1}^l z_j) - \sum_{j=1}^l F(z_j)| \leq C_{l,d} \sum_{j \neq j'} |z_j| |z'_j|^{4/d}$$
for some $C_{l,d} < \infty$ (see e.g. \cite[(1.10)]{bg}), and the claim follows from \eqref{ntg} (and \eqref{sbound}).
\end{proof}

Let $\delta > 0$ be a small number depending on $Bm_0$.  Then, by the above two lemmas and \eqref{sun} we have
$$ M( u^{(l)}_n(0) - u_n(0) ), \| (i \partial_t + \Delta) u^{(l)}_n - F(u^{(l)}_n) \|_{L^{2(d+2)/(d+4)}_{t,x}(\R \times \R^d)}
\leq \delta$$
and
$$ S(u^{(l)}_n) \leq 2 Bm_0, $$
provided $l$ is sufficiently large (depending on $\delta$) and $n$ sufficiently large (depending on $l,\delta$).  Applying Lemma \ref{stab} (for
$\delta$ chosen small enough depending on $2Bm_0$), we see that $u_n$ exists globally and
$$ S(u_n) \leq 3 Bm_0.$$
But this contradicts \eqref{under0}.

The only remaining possibility is that \eqref{mj} fails for every $\eps > 0$, and thus
$$  \sup_j M( \phi^{(j)} ) = m_0.$$
Comparing this with \eqref{un}, we see that at most one of the $\phi^{(j)}$ can be non-zero.  This means that the profile decomposition simplifies to
\begin{equation}\label{uno}
 u_n(0) = h_n e^{it_n \Delta} \phi + w_n
\end{equation}
for some sequence $t_n \in \R$ converging to $-\infty$, $0$, or $+\infty$, $h_n \in G$, some $\phi$ of mass $M(\phi) = m_0$, and some $w_n$ with
$M(w_n) \to 0$ (and hence $S(e^{it\Delta} w_n) \to 0$) as $n \to \infty$ (this is from \eqref{un}).  By applying the group action $T_{h_n^{-1}}$, we
may take $h_n$ to be the identity, and thus
$$ M( u_n(0) - e^{it_n \Delta} \phi ) \to 0 \hbox{ as } n \to \infty.$$

If $t_n$ converges to zero, then $u_n(0)$ is convergent in $L^2_x(\R^d)$ to $\phi$, and thus $G u_n(0)$ is convergent in $G \backslash L^2_x(\R^d)$,
as desired.  So the only remaining case is when $t_n$ converges to an infinite time; we shall assume that $t_n$ goes to $+\infty$, as the other case
is similar.  By \eqref{strich} we have
$$ S( e^{it\Delta} \phi ) < \infty$$
and hence, by time translation invariance and monotone convergence,
$$ \lim_{n \to \infty} S_{\geq 0}( e^{it\Delta} e^{it_n \Delta} \phi ) = 0.$$
Since the action of $G$ preserves linear solutions of the Schr\"odinger equation, we have $e^{it\Delta} h_n = T_{h_n} e^{it\Delta}$.  Since $T_{h_n}$
preserves the scattering norm $S$ (as well as $S_{\geq 0}$ and $S_{\leq 0}$) we thus have
$$ \lim_{n \to \infty} S_{\geq 0}( e^{it\Delta} h_n e^{it_n \Delta} \phi ) = 0.$$
Since $S(e^{it\Delta} w_n) \to 0$ as $n \to \infty$, we thus see from \eqref{uno} that
$$ \lim_{n \to \infty} S_{\geq 0}( e^{it\Delta} u_n(0) ) = 0.$$
Applying Lemma \ref{stab} (using $0$ as the approximate solution and $u_n(0)$ as the initial data) we conclude that
$$ \lim_{n \to \infty} S_{\geq 0}( u_n ) = 0.$$
But this contradicts one of the estimates in \eqref{under0}.  A similar argument allows us to exclude the case when $t_n$ goes to $-\infty$, where we
now use the other half of \eqref{under0}. This concludes the proof of Proposition \ref{mainprop}.

\section{A Duhamel formula}

We have just established Theorem \ref{main}, which shows that the critical mass $m_0$ is linked to the existence of maximal-lifespan solutions which
blow up both forward and backward in time and which are almost periodic modulo $G$.  It is thus of interest to study such solutions further.  We will
not do so in depth here (but see \cite{tvz}); however, we will note a Duhamel formula which may have application towards such a study.

Recall from Theorem \ref{lwp}(iii) that if a maximal-lifespan solution $u: I \times \R^d \to \C$ does not blow up forward in time, then $u$ scatters
forward in time to a linear solution $e^{it\Delta} u_+$, or equivalently that $e^{-it\Delta} u(t)$ is strongly convergent in $L^2_x(\R^d)$ as $t \to
+\infty$ to a scattering state $u_+$, which has the same mass as $u$.  Similar statements of course hold backward in time.  In sharp contrast, we
shall see that almost periodic (modulo $G$) solutions will exhibit the opposite behavior, namely that they are \emph{asymptotically orthogonal} to
all linear solutions at the endpoints of their lifespan $I$.  Thus, there is a dichotomy between scattering solutions and almost periodic solutions,
analogous to the distinction between free states and bound states in the study of linear Schr\"odinger equations with potential, except for the fact
that the almost periodicity is only relative to the symmetry group $G$.

More precisely, we have

\begin{proposition}\label{uconv}
Let $u: I \times \R^d \to \C$ be a maximal-lifespan solution which is almost periodic.  Then $e^{-it\Delta} u(t)$ is weakly convergent in
$L^2_x(\R^d)$ to zero as $t \to \sup(I)$ or $t \to \inf(I)$.
\end{proposition}

\begin{proof} Let us just prove the claim as $t \to \sup(I)$, as the reverse claim is similar.
By almost periodicity, we have a compact subset $K \subset L^2_x(G)$
and group elements $g_{\theta(t), \xi_0(t), x_0(t), \lambda(t)} \in G$ for each $t \in I$ such that
\begin{equation}\label{guk}
g_{\theta(t), \xi_0(t), x_0(t), \lambda(t)}^{-1} u(t) \in K.
\end{equation}
Suppose first that $\sup(I)$ is finite, and thus $u$ exhibits forward blowup in finite time.  We claim that this forces $\lambda(t)$ to go to zero as
$t \to \sup(I)$.  For if this were not the case, we could find a sequence $t_n \in I$ of times going to $\sup(I)$ such that $\lambda(t_n)$ is bounded
away from zero. But then, by the compactness of $K$ (and passing to a subsequence if necessary) we may ensure that $g_{\theta(t_n), \xi_0(t_n),
x_0(t_n), \lambda(t_n)}^{-1} u(t_n)$ is strongly convergent in $L^2_x(\R^d)$.  But as $t_n$ is converging to $\sup(I)$ and $\lambda(t_n)$ is bounded
away from zero, we see that the forward lifespan of solutions with this initial data converges to zero as $n \to \infty$.  This contradicts
Theorem~\ref{lwp}(i), (vi).

Since $\lambda(t) \to 0$ as $t \to \sup(I)$, the operators $g_{\theta(t), \xi_0(t), x_0(t), \lambda(t)}$ are weakly convergent to zero. By the
compactness of $K$, this implies
$$ \lim_{t \to \sup(I)} \sup_{f \in K} |\langle f, g_{\theta(t), \xi_0(t), x_0(t), \lambda(t)} \phi \rangle_{L^2_x(\R^d)}| = 0$$
for all $\phi \in L^2_x(\R^d)$.  From this and \eqref{guk}, we see that $u(t)$ converges weakly to zero as $t \to \sup(I)$.  Since $\sup(I)$ is
finite and the propagator curve $t \mapsto e^{it\Delta}$ is continuous in the strong operator topology, we see that $e^{-it\Delta} u(t)$ converges
weakly to zero, as desired.

Now suppose instead that $\sup(I)$ is infinite.  It will suffice to show that
$$ \lim_{t \to +\infty} \langle e^{-it\Delta} u(t), \phi \rangle_{L^2_x(\R^d)} = 0$$
for all test functions $\phi \in C^\infty_0(\R^d)$.  Applying \eqref{guk} and duality, it suffices to show that
$$ \lim_{t \to +\infty} \sup_{f \in K} |\langle f, g_{\theta(t), \xi_0(t), x_0(t), \lambda(t)} e^{it\Delta} \phi \rangle_{L^2_x(\R^d)}| = 0;$$
by the compactness of $K$, it therefore suffices to show that
$$ \lim_{t \to +\infty} \langle f, g_{\theta(t), \xi_0(t), x_0(t), \lambda(t)} e^{it\Delta} \phi \rangle_{L^2_x(\R^d)} = 0$$
for each $f \in L^2_x(\R^d)$ separately.  By density arguments we may take $f$ to also be a test function.  But the claim now follows from the
stationary phase expansion of $e^{it\Delta} \phi$ (or by using the fundamental solution), the point being that $e^{it\Delta} \phi$ acquires a
quadratic phase oscillation as $t \to \infty$ which cannot be renormalized by any of the symmetries in $G$.
\end{proof}

Let $u: I \times \R^d \to \C$ be as in the above proposition.  Recall the Duhamel formula
$$ u(t) = e^{it\Delta} e^{-it_+ \Delta} u(t_+) + i \int_t^{t_+} e^{i(t-t') \Delta} F( u(t') )\ dt'$$
for any $t, t_+ \in I$.  Letting $t_+$ converge to $\sup(I)$ and using Proposition~\ref{uconv}, we conclude the backward (advanced) Duhamel formula
$$ u(t) = i \int_t^{\sup(I)} e^{i(t-t') \Delta} F( u(t') )\ dt',$$
where the improper integral is interpreted in a conditionally convergent sense in the weak topology, that is,
$$ \langle u(t), f \rangle = \lim_{t_+ \to \sup(I)} \langle i \int_t^{t_+} e^{i(t-t') \Delta} F( u(t') )\ dt', f \rangle_{L^2_x(\R^d)}$$
for all $f \in L^2_x(\R^d)$.  Note that these partial integrals are equal to $u(t) - e^{it\Delta} e^{-it_+\Delta} u(t_+)$ and are thus uniformly
bounded in $L^2_x(\R^d)$.  Similarly, we have the forward (retarded) Duhamel formula
$$ u(t) = -i \int_{\inf(I)}^t e^{i(t-t') \Delta} F( u(t') )\ dt'.$$
One can interpolate between these two Duhamel formulae using an arbitrary operator $P: L^2_x(\R^d) \to L^2_x(\R^d)$ and obtain
a two-way (retarded-advanced) Duhamel formula
$$ u(t) = i P \int_t^{\sup(I)} e^{i(t-t') \Delta} F( u(t') )\ dt' - i (1-P) \int_{\inf(I)}^t e^{i(t-t') \Delta} F( u(t') )\ dt'.$$
For instance, in the spherically symmetric case, one might choose $P$ to be a pseudodifferential projection to ``outgoing'' waves, thus ensuring that
this formula expresses $u(t)$ in terms of integrals which mostly avoid the spatial origin $x=0$, which is presumably the most singular location for a
radial solution (cf. \cite{tao}). These formulae show that an almost periodic blowup solution is ``non-radiating'' or ``self-perpetuating''; the
solution can be expressed in terms of its own nonlinear interactions in the past and/or future, with no radiation terms coming from the endpoints of
the lifespan interval $I$.  Such formulae may be particularly useful in higher dimensions $d \geq 4$, when the decay of the fundamental solution (and
the subquadratic power of the nonlinearity $F$, coupled with mass conservation) ensures that these integrals are in fact uniformly convergent at the
endpoints $\inf(I)$, $\sup(I)$.  Variants of such formulae have appeared recently (see \cite{gopher}, \cite{tao:radialfocus}, \cite{tao}) and it is
quite possible that they can be used to provide further regularity and decay on almost periodic blowup solutions\footnote{These formulae are
analogous to the reproducing formula $Q = (-\Delta+1)^{-1} F(Q)$ for the ground state equation \eqref{ground}, which can be used iteratively to
establish arbitrary amounts of smoothness and decay for this ground state.}.

\section{Spherically symmetric analogues}\label{spherical-sec}

In this section, we specialize the study of the equation \eqref{nls} to spherically symmetric initial data $u_0(x) = u_0(|x|)$, and hence (by
rotation symmetry and uniqueness) to spherically symmetric solutions $u(t,x) = u(t,|x|)$.  Of course, one expects the theory here to be significantly
simpler than in the general case, though we will also reduce the size of the symmetry group $G$ and so the deduction of the results in this section
from the preceding material is not entirely trivial.

Let $L^2_x(\R^d)_\rad$ be the closed subspace of $L^2_x(\R^d)$ consisting of spherically symmetric functions.  The well-posedness theory in Theorem
\ref{lwp} has an obvious counterpart (which we will not restate here), in which all functions are spherically symmetric.

The full group $G$ no longer acts on $L^2_x(\R^d)_\rad$, and so we shall work instead with the maximal subgroup of $G$ which does preserve
$L^2_x(\R^d)_\rad$, namely the phase rotations\footnote{Actually the phase rotations now serve no useful role and can be discarded here if desired.
Even in the non-symmetric case, the phase rotations were only necessary because they arose from commutators of the spatial translations and frequency
modulations, but did not actually supply any new sources of non-compactness.} and dilations:
$$G_\rad := \{ g_{\theta,0,0,\lambda}: \theta \in \R/2\pi\Z, \lambda > 0\}.$$
Similarly, we have the spherically symmetric enlarged group
$$G'_\rad := \{ g_{\theta,0,0,\lambda,t_0}: t_0 \in \R, \theta \in \R/2\pi\Z, \lambda > 0\}.$$
We define the radial analogue $A(m)_\rad$ of $A(m)$ by using \eqref{Adef} as before, but now restricting $u$ to spherically symmetric solutions;
thus, $A(m)_\rad \leq A(m)$.  We can then define the \emph{spherically symmetric critical mass} $0 < m_{0,\rad} \leq +\infty$ to be the unique value
such that $A(m)_\rad$ is finite for $m < m_{0,\rad}$ and infinite for $m \geq m_{0,\rad}$; hence, $m_{0,\rad} \geq m_0$.  Since the ground state $Q$
is already spherically symmetric, we know that $m_{0,\rad} \leq M(Q)$ in the focusing case $\mu=+1$.  Thus, Conjecture~\ref{conj} has a spherically
symmetric counterpart:

\begin{conjecture}[Scattering conjecture, spherically symmetric case]\label{conj-sph}
In the defocusing case ($\mu=+1$) we have $m_{0,\rad}=+\infty$, while in the focusing case ($\mu=-1$) we have $m_{0,\rad} = M(Q)$.
\end{conjecture}

In the sequel to this paper, \cite{tvz}, we shall verify this weaker conjecture in high dimensions $d \geq 3$ and with defocusing sign $\mu=+1$.

We can define the concept of being an almost periodic solution modulo $G_\rad$ in obvious analogy with being almost periodic modulo $G$.  The main
result is then:

\begin{theorem}[Reduction to almost periodic solutions, spherically symmetric case]\label{main-rad}
Fix $\mu$ and $d$ and suppose that the spherically symmetric critical mass $m_{0,\rad}$ is finite.  Then there exists a maximal-lifespan spherically
symmetric solution $u$ of mass exactly $m_{0,\rad}$ which blows up both forward and backward in time.  Furthermore, any maximal-lifespan solution $u:
I \times \R^d \to \C$ which blows up both forward and backward in time is almost periodic modulo $G_\rad$, and there exists a function $N: I \to
(0,+\infty)$ such that for every $\eta > 0$ there exists $0 < C(\eta) < \infty$ (possibly depending on $u$) such that
\begin{equation}\label{spaceconc-rad}
\int_{|x| \geq C(\eta)/N(t)} |u(t,x)|^2\ dx \leq \eta
\end{equation}
and
\begin{equation}\label{freqconc-rad}
\int_{|\xi| \geq C(\eta) N(t)} |\hat u(t,\xi)|^2\ d\xi \leq \eta
\end{equation}
for all $t \in I$.
\end{theorem}

By repeating the arguments in the previous sections, restricting all functions to be spherically symmetric, we see that to prove this theorem it
suffices to establish the spherically symmetric counterpart of Theorem~\ref{lp}, namely

\begin{theorem}[Spherically symmetric linear profiles]\label{lp-sym}
Fix $d$. Let $u_n$, $n=1,2,\ldots$ be a bounded sequence in $L^2_x(\R^d)_\rad$.  Then (after passing to a subsequence if necessary) there exists a
family $\phi^{(j)}$, $j=1,2,\ldots$ of functions in $L^2_x(\R^d)_\rad$ and group elements $g^{(j)}_n \in G'_\rad$ for $j,n = 1,2,\ldots$ such that we
have the decomposition \eqref{undecomp} for all $l=1,2,\ldots$, where $w^{(l)}_n \in L^2_x(\R^d)_\rad$ obeys \eqref{sln}.  Furthermore, $g^{(j)}_n,
g^{(j')}_n$ are asymptotically orthogonal for any $j \neq j'$, and for any $l \geq 1$ we have the mass decoupling property \eqref{un}.
\end{theorem}

\begin{proof}
We apply Theorem \ref{lp} to obtain a preliminary decomposition with most of the desired properties.  The ones which are missing are that $g^{(j)}_n$
lie in the large group $G'$ rather than in the smaller group $G'_\rad$, and that the functions $\phi^{(j)}$ and $w^{(l)}_n$ are not spherically
symmetric.  We will thus need to perturb these objects slightly to obtain the desired symmetry.

Let $P: L^2_x(\R^d) \to L^2_x(\R^d)$ be the orthogonal projection onto $L^2_x(\R^d)_\rad$; note that this operator commutes with the linear
propagators $e^{it\Delta}$ and with the group $G'_\rad$.  Applying $1-P$ to \eqref{undecomp} and using the spherical symmetry of $u_n$ we obtain for
any $l=1,2,\ldots$
\begin{equation}\label{oo}
 0 = \sum_{j=1}^l (1-P) g^{(j)}_n \phi^{(j)} + (1-P) w^{(l)}_n.
 \end{equation}
From the asymptotic orthogonality of the $g^{(j)}_n$ we already have
$$ \lim_{n \to \infty} \langle g^{(j)}_n \phi^{(j)}, g^{(j')}_n \phi^{(j')} \rangle_{L^2_x(\R^d)} = 0$$
for $j \neq j'$.  A variant of the same argument also gives
\begin{equation}\label{ao1}
\lim_{n \to \infty} \langle (1-P) g^{(j)}_n \phi^{(j)}, (1-P) g^{(j')}_n \phi^{(j')} \rangle_{L^2_x(\R^d)} = 0;
\end{equation}
this is easiest to establish by first considering the case when $\phi^{(j)}, \phi^{(j')} \in C^\infty_0(\R^d)$ are test functions, and then using
density to obtain the general case.  Similarly, from \eqref{sln} and the symmetries of the linear Schr\"odinger equation we have
$$ \lim_{n \to \infty} S( e^{it\Delta} (g^{(j)}_n)^{-1} (1-P) w^{(l)}_n ) = 0$$
for any $j,l$, and hence $(g^{(j)}_n)^{-1} (1-P) w^{(l)}_n$ converges weakly to zero.  In particular,
\begin{equation}\label{ao2}
\lim_{n \to \infty} \langle (1-P) g^{(j)}_n \phi^{(j)}, (1-P) w^{(l)}_n \rangle_{L^2_x(\R^d)} = 0.
\end{equation}
Comparing \eqref{ao1} and \eqref{ao2} against \eqref{oo}, we conclude
\begin{equation}\label{pj}
\lim_{n \to \infty} M( (1-P) g^{(j)}_n \phi^{(j)} ) = 0
\end{equation}
for all $j$, and
$$ \lim_{n \to \infty} M( (1-P) w^{(l)}_n ) = 0$$
for all $l$.  Thus $g^{(j)}_n \phi^{(j)}$ and $w^{(l)}_n$ are asymptotically spherically symmetric.

Suppose that $j$ is such that $\phi^{(j)}$ is non-zero.  We now claim that $g^{(j)}_n$ stays close to $G'_\rad$ in the following sense: if
$$ g^{(j)}_n = g_{\theta^{(j)}_n, x^{(j)}_n, \xi^{(j)}_n, \lambda^{(j)}_n, t^{(j)}_n },$$
then we claim that $x^{(j)}_n /\lambda^{(j)}_n$ and $\xi^{(j)}_n  \lambda^{(j)}_n$ are bounded in $n$. To see this, suppose for contradiction that
this were not the case; then, by passing to a subsequence if necessary we can find an $j$ such that
$$ \lim_{n \to \infty} |x^{(j)}_n|/ \lambda^{(j)}_n + |\xi^{(j)}_n|  \lambda^{(j)}_n = +\infty.$$
One can easily verify that this forces $g^{(j)}_n \phi^{(j)}$ to be asymptotically non-radial, and more precisely that
$$ \lim_{n \to \infty} M( P g^{(j)}_n \phi^{(j)} ) = 0$$
(indeed, one can verify this first for test functions $\phi^{(j)} \in C^\infty_0(\R^d)$ and then use density).  Since we are assuming $\phi^{(j)}$ to
be non-zero, we contradict \eqref{pj}.  Thus the sequences $x^{(j)}_n \lambda^{(j)}_n$ and $\xi^{(j)}_n / \lambda^{(j)}_n$ are bounded.  By passing
to a subsequence for each $j$ and using the usual diagonalisation argument, we may in fact assume that $x^{(j)}_n \lambda^{(j)}_n$ and $\xi^{(j)}_n /
\lambda^{(j)}_n$ converge for each $j$.  By applying an appropriate element of $G$ to $\phi^{(j)}$, we may in fact assume that these quantities
converge to zero.  One can then \emph{replace} these quantities with zero, absorbing the errors into $w^{(l)}_n$ (which will be acceptable by
\eqref{strich}).  To summarize, we can now assume that $x^{(j)}_n = \xi^{(j)}_n = 0$, or in other words that $g^{(j)}_n \in G'_\rad$.  In this case
$g^{(j)}_n$ commutes with $P$, and hence \eqref{pj} simplifies to
$$ M( (1-P) \phi^{(j)} ) = 0.$$
Thus, $\phi^{(j)}$ is spherically symmetric.

In the degenerate cases when $\phi^{(j)} = 0$, we of course already have $\phi^{(j)}$ spherically symmetric, and we can easily also ensure
$g^{(j)}_n$ to lie in $G'_\rad$ without disrupting the pairwise asymptotic orthogonality.  From \eqref{undecomp} we now conclude that $w^{(l)}_n$ is
also spherically symmetric, and we are done.
\end{proof}

\begin{remark}
One could also prove Theorem \ref{lp-sym} by modifying the arguments in \cite{begout}, \cite{mv} directly.  Indeed, this would be the most direct
approach (for instance, the bilinear restriction theorem from \cite{tao:bilinear} is significantly easier to prove if one assumes spherical
symmetry, and one can use local smoothing estimates and radial Sobolev inequalities, such as those in \cite{vilela}, to significantly increase the
amount of compactness available).  However, this would require a large amount of tedious repetition of existing arguments in the literature and so we
will not do so here.
\end{remark}

\section{Frequency and space localization -- a quantitative version}\label{quantitative}

In this section we outline the quantitative proof of frequency and space localization.  The argument follows closely the ones in
\cite{bourg.critical, gopher, rv, visan:art}.  We will need a few small parameters,
$$
1 \gg \eta_0 \gg \eta_1 \gg \eta_2 > 0
$$
where each $\eta_j$ is allowed to depend on the critical mass $m_0$ and on any of the larger $\eta$'s. We will choose
$\eta_j$ small enough such that, in particular, it will be smaller than any constant depending on the
previous $\eta$'s used in the argument.

Rather than working with a minimal-mass blowup solution, in this section we will work with a minimal-mass \emph{almost} blowup solution.
More precisely, fix $\mu$ and $d$ and assume that the critical mass $m_0$ is finite.  Then, we make the following
\begin{definition}\label{almost}
A \emph{minimal-mass almost blowup solution} to \eqref{nls} is a solution $u$ on a compact time interval $I_*$ such that
$$
M(u)= m_0
$$
and
\begin{align}\label{HUGE}
S(u) >1/\eta_2.
\end{align}
\end{definition}

In the spirit of \cite{bourg.critical, gopher, rv, visan:art}, the quantitative proof of localization relies heavily on induction on mass techniques.
If a minimal-mass almost blowup solution were not localized in both physical and frequency space, it could be decomposed into two essentially
separate solutions, each with strictly smaller mass than the original.  As $A(m)<\infty$ for $m<m_0$, we can then extend these smaller mass
solutions to all of $I_*$.  As each of the separate evolutions exactly solves \eqref{nls}, we expect their sum to solve \eqref{nls}
approximately.  We could then use perturbation theory to derive a bound on $S(u)$, thus contradicting the fact that $\eta_2$ can be chosen arbitrarily
small in \eqref{HUGE}.

In the remainder of this section, we will present this argument in more detail.  We begin with frequency localization.

\begin{proposition}[Frequency delocalization implies spacetime bounds]\label{no delocal}
Let $\eta>0$ and suppose there exists a dyadic frequency $N_0>0$ and a time $t_0\in I_*$ such that we have the mass separation conditions
$$
\|P_{\le N_0}u(t_0)\|_{L_x^2}\ge \eta,
$$
$$
\|P_{\ge K(\eta)N_0}u(t_0)\|_{L_x^2}\ge \eta.
$$
Then, if $K(\eta)$ is sufficiently large depending on $\eta$, we have
$$
\|u\|_{L_{t,x}^{\frac{2(d+2)}{d}}(I_*\times\R^d)}\le C(\eta).
$$
\end{proposition}

Using this proposition and Definition~\ref{almost} (taking $\eta_2$ sufficiently small in \eqref{HUGE}), we immediately see that a minimal-mass
almost blowup solution must be localized in frequency (see \cite{gopher, rv, visan:art} for the analogue statement in the energy-critical setting).
\begin{corollary}[Frequency localization]\label{freq local}
Let $u$ be a minimal-mass almost blowup solution on $I_*\times\R^d$.  Then, for each time $t\in I_*$, there exists a dyadic frequency
$N(t)\in 2^Z$ such that for every $\eta_1\le \eta\le \eta_0$, we have small mass at frequencies $\ll N(t)$
$$\|P_{\le c(\eta)N(t)}u(t)\|_{L_x^2}\le \eta,$$
small mass at frequencies $\gg N(t)$
$$
\|P_{\ge C(\eta)N(t)}u(t)\|_{L_x^2}\le \eta,
$$
and large mass at frequencies $\sim N(t)$.
$$
\|P_{c(\eta)N(t)<\cdot<C(\eta)N(t)}u(t)\|_{L_x^2}\sim 1.
$$
Here, the values $0<c(\eta)\ll 1 \ll C(\eta)$ depend on $\eta$.
\end{corollary}

\emph{Sketch of proof of Proposition~\ref{no delocal}}.
The proof of this proposition follows the same strategy used to derive its analogue in the energy-critical setting.
Using mass conservation and the pigeonhole principle, we can find a frequency band where the solution has very little mass.
Taking $K(\eta)$ very large and rescaling appropriately, we may assume that
$$
\|P_{\eps^2\leq \cdot \leq \eps^{-2}} u(t_0)\|_{L_x^2}\leq \eps,
$$
where $0<\eps=\eps(\eta)\ll 1$ will be chosen later.  We then define $u_{lo}(t_0):=P_{\leq \eps}u(t_0)$ and $u_{hi}(t_0):=P_{\geq \eps^{-1}}u(t_0)$.
By hypothesis, we immediately see that
\begin{align*}
M(u_{lo}(t_0)), M(u_{hi}(t_0)) \leq m_0-\eta^2,
\end{align*}
and thus we can find two global solution $u_{lo}$, $u_{hi}$ to \eqref{nls} with initial data $u_{lo}(t_0)$, $u_{hi}(t_0)$ respectively, such that
\begin{align*}
S(u_{lo}), S(u_{hi})\leq A(m_0-\eta^2)=C(\eta).
\end{align*}
By Lemma~\ref{stab}, we thus see that the claim of Proposition~\ref{no delocal} would follow immediately (taking $\eps=\eps(\eta)$ sufficiently small)
if we could show that $\tilde u:= u_{lo}+u_{hi}$ is an approximate solution to \eqref{nls} in the sense that
\begin{align}\label{approx sol}
\|i\tilde u_t +\Delta \tilde u - F(\tilde u)\|_{L_{t,x}^{\frac{2(d+2)}{d+4}}(I_*\times\R^d)}\leq C(\eta)\eps^c,
\end{align}
for some constant $c>0$.  In order to prove \eqref{approx sol}, one has to control interactions between $u_{lo}$ and $u_{hi}$.
This is done by showing that $u_{hi}$ remains essentially at high frequencies, $u_{lo}$ remains essentially at low frequencies, and using this
information to control the interactions.  In the energy-critical setting, one uses the conservation of mass, the persistence of positive regularity,
and the bilinear Strichartz estimate respectively (see \cite{gopher, rv, visan:art}).
In the mass-critical setting, these tools are replaced by persistence of negative regularity (see Lemma~\ref{keep negative} below), persistence of
positive regularity, and the bilinear restriction estimate below.

\begin{lemma}[Bilinear restriction estimate,  \cite{tao:bilinear}]\label{q bilinear}
Let $I$ be a compact time interval, $t_0 \in I$, $N>0$, and let $u_1,\ u_2$ be two solutions to \eqref{nls} such that $u_j(t)$ has Fourier transform
supported in the region $ \{|\xi_j|\le N\}$ for $j=1,2$.  Suppose also that the Fourier supports of $u_1,\ u_2$ are separated by at least $\ge cN$.
Then, for any $q>\frac{d+3}{d+1}$ we have
$$
\|u_1u_2\|_{L_{t,x}^q(I\times\R^d)}\lesssim  N^{d-\frac{d+2}q}\|u_1\|_{S_*^0(I\times\R^d)}\|u_2\|_{S_*^0(I\times\R^d)},
$$
where $S_*^0$ is the strong Strichartz norm
$$
\|u\|_{S_*^0(I\times\R^d)}=\|u(t_0)\|_2+\|(i\partial_t+\Delta)u\|_{L_{t,x}^{\frac{2(d+2)}{d+4}}(I\times\R^d)}.
$$
\end{lemma}

We will also rely on this bilinear restriction estimate to prove

\begin{lemma}[Persistence of negative regularity]\label{keep negative}
Fix $\mu$ and $d$ and let $u$ be a solution to \eqref{nls} on a time interval  $I=[t_0,t_1]$ such that
$$
\|u(t_0)\|_{L_x^2}\leq M,
$$
and
\begin{align}\label{L-bound}
\|u\|_{L_{t,x}^{\frac{2(d+2)}{d}}(I\times\R^d)}\le L.
\end{align}
Suppose also that for all $N\in 2^{\mathbb Z}$, some constant $A>0$, and some small constant $s>0$ we have
\begin{align}\label{nr-hyp}
\|P_Nu(t_0)\|_{L_x^2} \le AN^s.
\end{align}
Then
\begin{align}\label{nr-claim}
\|P_N u\|_{L_{t,x}^{\frac{2(d+2)}{d}}(I\times\R^d)}\leq C(L,M)AN^s.
\end{align}
\end{lemma}

\begin{proof}
Subdividing the time interval $I$, we see that we need only prove Lemma~\ref{keep negative} with hypothesis \eqref{L-bound} being replaced by
\begin{align}\label{eta-bound}
\|u\|_{L_{t,x}^{\frac{2(d+2)}{d}}(I\times\R^d)}\le \eta,
\end{align}
for a small constant $\eta>0$ to be chosen later.  Throughout the rest of the proof, all spacetime norms will be on $I\times\R^d$.

Taking $\eta=\eta(M)$ sufficiently small, by Strichartz we estimate
\begin{align}\label{us0*}
\|u\|_{S^0_*}
&\lesssim \|u(t_0)\|_{L_x^2} + \||u|^{\frac 4d}u\|_{L_{t,x}^{\frac{2(d+2)}{d+4}}}
\lesssim M + \|u\|_{L_{t,x}^{\frac{2(d+2)}{d}}}^{1+\frac 4d}
\lesssim M.
\end{align}

To establish \eqref{nr-claim}, it suffices to prove
\begin{align}\label{nr-claim-eps}
\|P_N u\|_{L_{t,x}^{\frac{2(d+2)}{d}}}\leq C(M)(AN^s+\eps), \quad \text{for all }  N>0
\end{align}
for any $\eps>0$ since shrinking $\eps$ to zero yields the claim.  By standard continuity arguments, it suffices to prove \eqref{nr-claim-eps}
under the additional assumption
\begin{align}\label{nr-ass}
\|P_N u\|_{L_{t,x}^{\frac{2(d+2)}{d}}}\leq C_0(AN^s+\eps), \quad \text{for all }  N>0
\end{align}
for some large constant $C_0=C_0(M)>0$.

By Strichartz and \eqref{nr-hyp}, it suffices to prove that
\begin{equation*}
\|P_N\int_{t_0}^t e^{i(t-s)\Delta} (|u|^{\frac 4d}u)(s)ds\|_{L_{t,x}^{\frac{2(d+2)}d}}\lesssim_M AN^s + \eps.
\end{equation*}
By duality, we reduce to showing that
\begin{align}\label{keep negative1}
\|G_N |u|^{\frac 4d}u\|_{L_{t,x}^1}\lesssim_M AN^s + \eps
\end{align}
where $G:I\times\R^d \to \C$ ranges over all functions such that
$$
\|G\|_{L_{t,x}^{\frac{2(d+2)}{d+4}}}\leq 1
$$
and
$$
G_N(x,t):=\int_t^{t_1}P_N e^{i(t-s)\Delta}G(x,s)\,ds.
$$
From Strichartz, we immediately see that
\begin{align}\label{sgn}
\|G_N\|_{S^0_*}\lesssim \|G\|_{L_{t,x}^{\frac{2(d+2)}{d+4}}}\lesssim 1.
\end{align}

Decomposing $u:=u_{\lesssim N} + u_{\gg N}$ and using the triangle inequality, we see that \eqref{keep negative1} would follow from
\begin{align}
\|G_N |u_{\lesssim N}|^{1+\frac 4d}\|_{L_{t,x}^1}&\lesssim_M AN^s + \eps \label{ts-lo}\\
\|G_N |u_{\gg N}|^{1+\frac 4d}\|_{L_{t,x}^1}&\lesssim_M AN^s + \eps. \label{ts-hi}
\end{align}

To prove \eqref{ts-lo}, we use H\"older, \eqref{eta-bound}, \eqref{nr-ass}, \eqref{sgn}, and take $\eta$ sufficiently small:
\begin{align*}
\|G_N |u_{\lesssim N}|^{1+\frac 4d}\|_{L_{t,x}^1}
\lesssim \|G\|_{L_{t,x}^{\frac{2(d+2)}{d+4}}}\|u_{\lesssim N}\|_{L_{t,x}^{\frac{2(d+2)}{d}}}^{1+\frac 4d}
\lesssim \eta^{\frac 4d} C_0 (A N^s + \eps)
\lesssim A N^s + \eps.
\end{align*}

To prove \eqref{ts-hi}, we will consider two cases: $d\leq 3$ and $d\geq 4$.  Let first $d\leq 3$ and let $p,q$ be such that
$\frac{d+3}{d+1}<p<\frac{d+2}d$ and
$$
\frac 1p +\frac 1q= \frac{3d}{2(d+2)}.
$$
Then, applying a dyadic decomposition to $u_{\gg N}$ and using H\"older's inequality and \eqref{eta-bound}, we get
\begin{align*}
\|G_N & |u_{\gg N}|^{1+\frac 4d}\|_{L_{t,x}^1}\\
&\lesssim \sum_{N_1\geq N_2 \gg N} \|G_N u_{N_1} u_{N_2} |u_{\gg N}|^{\frac 4d -1}\|_{L_{t,x}^1}\\
&\lesssim \sum_{N_1\geq N_2 \gg N} \|u_{N_1}\|_{L_{t,x}^{\frac{2(d+2)}d}}^{\frac 12} \|u_{\gg N}\|_{L_{t,x}^{\frac{2(d+2)}d}}^{\frac 4d -1} \|G_N u_{N_1}\|_{L_{t,x}^p}^{\frac 12} \|u_{N_2}\|_{L_{t,x}^{\frac{2(d+2)}d}} \|G_N\|_{L_{t,x}^q}^{\frac 12}\\
&\lesssim \eta^{\frac 4d -\frac 12}\sum_{N_1\geq N_2 \gg N} \|G_N u_{N_1}\|_{L_{t,x}^p}^{\frac 12} \|u_{N_2}\|_{L_{t,x}^{\frac{2(d+2)}d}} \|G_N\|_{L_{t,x}^q}^{\frac 12}.
\end{align*}
As for any $N_1\gg N$, the Fourier supports of $G_N$ and $u_{N_1}$ are separated by at least (say) $N_1/8$, we apply Lemma~\ref{q bilinear} followed
by \eqref{us0*} and \eqref{sgn} to get
\begin{align}
\|G_N u_{N_1}\|_{L_{t,x}^p}
\lesssim N_1^{d-\frac{d+2}p} \|G_N\|_{S^0_*}\|u_{N_1}\|_{S_*^0}
\lesssim_M N_1^{d-\frac{d+2}p}.\label{pgn}
\end{align}
On the other hand, by Sobolev embedding, Bernstein, and \eqref{sgn}, we have
\begin{align}
\|G_N\|_{L_{t,x}^q}
\lesssim \||\nabla|^{\frac d2 - \frac{d+2}q} G_N\|_{L_t^q L_x^{\frac{2dq}{dq-4}}}
\lesssim N^{\frac d2 - \frac{d+2}q} \|G_N\|_{S^0_*}
\lesssim N^{\frac d2 - \frac{d+2}q}. \label{qgn}
\end{align}
Thus,
\begin{align*}
\|G_N  |u_{\gg N}|^{1+\frac 4d}\|_{L_{t,x}^1}
\lesssim_M \eta^{\frac 4d -\frac 12}\sum_{N_1\geq N_2 \gg N} N_1^{\frac 12 (d-\frac{d+2}p)} C_0 (A N_2^s + \eps) N^{\frac 12(\frac d2 - \frac{d+2}q)}.
\end{align*}
Choosing $s$ sufficiently small such that
$$
\frac 12 (d-\frac{d+2}p)+s<0,
$$
summing first in $N_1$ and then in $N_2$, and choosing $\eta$ sufficiently small, we derive \eqref{ts-hi}.

We consider next the case $d\geq 4$.  Using again a dyadic decomposition, we estimate
\begin{align*}
\|G_N  |u_{\gg N}|^{1+\frac 4d}\|_{L_{t,x}^1}
&\lesssim \sum_{N_1\geq N_2 \gg N} \|G_N  |u_{N_1}|^{\frac 4d} u_{N_2}\|_{L_{t,x}^1} \\
&\quad + \sum_{N_1\geq N_2 \gg N} \|G_N u_{N_1} |u_{N_2}|^{\frac 4d} \|_{L_{t,x}^1}.
\end{align*}
Choosing $p,q$ such that $\frac{d+3}{d+1}<p<\frac{d+2}d$ and
$$
\frac 2p +\frac {d-2}q= \frac d2
$$
and using H\"older, \eqref{eta-bound}, \eqref{nr-ass}, \eqref{pgn}, \eqref{qgn}, and taking $\eta$ sufficiently small, we estimate
\begin{align*}
\sum_{N_1\geq N_2 \gg N} \|G_N & |u_{N_1}|^{\frac 4d} u_{N_2}\|_{L_{t,x}^1}\\
&\lesssim \sum_{N_1\geq N_2 \gg N} \|u_{N_1}\|^{\frac 2d}_{L_{t,x}^{\frac{2(d+2)}d}} \|G_N u_{N_1}\|^{\frac 2d}_{L_{t,x}^p} \|u_{N_2}\|_{L_{t,x}^{\frac{2(d+2)}d}} \|G_N\|_{L_{t,x}^q}^{1-\frac 2d}\\
&\lesssim_M \eta^{\frac 2d} \sum_{N_1\geq N_2 \gg N} N_1^{\frac 2d (d-\frac{d+2}p)} C_0 (A N_2^s + \eps) N^{(1-\frac 2d)(\frac d2 - \frac{d+2}q)}\\
&\lesssim_M A N^s + \eps.
\end{align*}
Similarly,
\begin{align*}
&\sum_{N_1\geq N_2 \gg N}  \| G_N u_{N_1}  |u_{N_2}|^{\frac 4d} \|_{L_{t,x}^1}\\
&\qquad \lesssim \sum_{N_1\geq N_2 \gg N} \|u_{N_1}\|^{\frac 2d}_{L_{t,x}^{\frac{2(d+2)}d}} \|G_N u_{N_1}\|^{\frac 2d}_{L_{t,x}^p} \|u_{N_1}\|^{1-\frac 4d}_{L_{t,x}^{\frac{2(d+2)}d}} \|u_{N_2}\|^{\frac 4d}_{L_{t,x}^{\frac{2(d+2)}d}} \|G_N\|_{L_{t,x}^q}^{1-\frac 2d}\\
&\qquad \lesssim_M \eta^{\frac 2d} \sum_{N_1\geq N_2 \gg N} N_1^{\frac 2d (d-\frac{d+2}p)} C_0 (A N_1^s + \eps)^{(1-\frac 4d)} (A N_2^s + \eps)^{\frac 4d} N^{(1-\frac 2d)(\frac d2 - \frac{d+2}q)}\\
&\qquad \lesssim_M A N^s + \eps.
\end{align*}
Thus, \eqref{ts-hi} holds for $d\geq 4$.  This concludes the proof of Lemma~\ref{keep negative}.
\end{proof}

Using the persistence of positive and negative regularity as well as Lemma~\ref{q bilinear}, one can control the interactions between $u_{lo}$
and $u_{hi}$ and prove that $\tilde u$ is an approximate solution in the sense of \eqref{approx sol}; we omit the details.
Taking $\eps=\eps(\eta)$ sufficiently small and using Lemma~\ref{stab}, we derive
$$
\|u\|_{L_{t,x}^{\frac{2(d+2)}{d}}(I_*\times\R^d)}\le C(\eta),
$$
thus concluding the proof of Proposition~\ref{no delocal}.

We turn next to space localization; the approach is that of Bourgain, \cite{bourg.critical}.
We split the interval $I_*$ into three consecutive subintervals ($I_{-}, I_0, I_{+}$) such that each of these subintervals carries a third of the total
$L_{t,x}^{\frac{2(d+2)}{d}}(I_*\times\R^d)$ mass of $u$.  It is on the middle subinterval $I_0$ that space localization is proved.
The proof is carried out in two steps: Firstly, one establishes space concentration, which basically means that the solution is big somewhere;
to prove space localization, one then has to show that the solution is small everywhere else.

In order to prove space concentration in the energy-critical setting, one first proves a lower bound on the
potential energy (which is a scale-invariant norm) on $I_0$ (see \cite{gopher, rv, visan:art}).  For the mass-critical NLS, we do something similar.  More precisely, we show
that for every $t\in I_0$ there exists $\xi(t)$ such that
\begin{align}\label{lb-gu}
\||\nabla|^{-\frac 14} G_{\xi(t)} (u) (t)\|_{L^{\frac{4d}{2d-1}}_x} \geq \eta_1,
\end{align}
where $G_{\xi}$ denotes the Galilean transform given by
$$
G_{\xi}(u)(t,x)=(T_{g_{0,\xi,0,1}}u)(t,x)=e^{ix\xi}e^{-it|\xi|^2}u(t,x-2t\xi).
$$
\begin{remark}
The norm on the left-hand side of \eqref{lb-gu} was chosen to be critical with respect to the scaling and dominated by the mass
(by Sobolev embedding); it is not the only norm with these properties and the exact choice is not essential for our argument.
\end{remark}

\begin{remark}
As Galilean transformations leave the equation \eqref{nls}, the mass, and the Strichartz norm $S$ invariant, we see that the Galilean transform
of a minimal-mass almost blowup solution is still a minimal-mass almost blowup solution.  In particular, it is still frequency-localized.
Moreover, because Galilean transformations leave the equation \eqref{nls} invariant, but not norms such as the one on the left-hand side
of \eqref{lb-gu}, it is not reasonable to expect \eqref{lb-gu} to hold without the presence of a Galilean transformation.
\end{remark}

The argument used to establish \eqref{lb-gu} is inspired by its analogue in \cite{gopher, rv, visan:art} and we outline it next.
The proof is by contradiction.  We assume there exists $t_0 \in I_0$ such that for any Galilean transformation $G_{\xi}$ we have
\begin{align}\label{no-lb-gu}
\||\nabla|^{-\frac 14} G_{\xi} (u) (t_0)\|_{L^{\frac{4d}{2d-1}}_x} < \eta_1.
\end{align}
We first note that the $L_{t,x}^{\frac{2(d+2)}{d}}$ norm of the free evolution of $u(t_0)$ must be large, since otherwise perturbation
theory would imply a bound on the scattering size $S(u)$ and choosing $\eta_2$ sufficiently small in \eqref{HUGE} would yield a contradiction.
Now, by Corollary~\ref{freq local}, the $L_{t,x}^{\frac{2(d+2)}{d}}$ norm of the free evolution of $u(t_0)$ can only be large on an annulus
centered around frequency $N(t_0)$.  Using an inverse Strichartz theorem (see the appendix in \cite{tao-lens}), we deduce that there exist
$t_1\in \R$, $x_1, \xi_1\in \R^d$, and $N_1\in 2^{\Z}$ such that
\begin{align*}
\int_{|\xi-\xi_1|\le N_1C(\eta_0)}|\widehat{P_{med}u(t_0)}(\xi)|^2\, d\xi &\gtrsim \eta_0\\
\int_{|x-x_1|\leq \frac{C(\eta_0)}{N_1}}|e^{i(t_1-t_0)\Delta} P_{med} u(t_0,x)|^2\, dx &\gtrsim \eta_0,
\end{align*}
where $P_{med}:=P_{c(\eta_0)N(t_0)<\cdot <C(\eta_0)N(t_0)}$.  In other words, the free evolution of $P_{med} u(t_0)$ can only be large in
$L_{t,x}^{\frac{2(d+2)}{d}}$ if it concentrates at some point.

Next, we use a Galilean transformation to send $\xi_1$ to zero.  More precisely, we define $\xi(t_0):=-\xi_1$ and let $\tilde u:= G_{\xi(t_0)}(u)$.
Let $\tilde N (t_0)$ be the frequency at which the new minimal-mass almost blowup solution $\tilde u$ is localized at time $t_0$, and denote
$\tilde P_{med}:=P_{c(\eta_0)\tilde N(t_0)<\cdot <C(\eta_0)\tilde N(t_0)}$.
Trivial computations yield
\begin{align}
\int_{|\xi|\le N_1C(\eta_0)}|\widehat{\tilde P_{med}\tilde u(t_0)}(\xi)|^2\, d\xi &\gtrsim \eta_0 \notag \\
\int_{|x-x_1 + 2t_1\xi_1|\leq \frac{C(\eta_0)}{N_1}}|e^{i(t_1-t_0)\Delta} \tilde P_{med} \tilde u(t_0,x)|^2\, dx &\gtrsim \eta_0.\label{smthg}
\end{align}
In particular, this implies that $c(\eta_0)\tilde N(t_0)<N_1 <C(\eta_0)\tilde N(t_0)$.  Given the dispersive effect of the free Schr\"odinger evolution,
\eqref{no-lb-gu} and \eqref{smthg} imply that $t_0$ and $t_1$ must be far apart (see, for example, Section~6 in \cite{gopher}).

We proceed next to remove the linear evolution of this bubble of concentration.  More precisely, let
$f(t):= e^{i(t-t_1)\Delta}f(t_1)$ with
$$
f(t_1):=\tilde P_{med}\bigl(\chi_{|x-x_1 + 2t_1\xi_1|\leq \frac{C(\eta_0)}{N_1}}e^{i(t_1-t_0)\Delta}\tilde P_{med}\tilde u(t_0)\bigr)
$$
and write $\tilde u(t_0)=v(t_0) + \alpha f(t_0)$ where $\alpha\in\C$ and $\langle v(t_0),f\rangle=0$.
By \eqref{smthg}, we have
$$
\langle \tilde u(t_0),f(t_0)\rangle \geq \eta_0,
$$
and so
$$
M(v(t_0))\leq M(\tilde u(t_0)) - c(\eta_0) < m_0.
$$
Hence, there exists a unique global solution $v$ to \eqref{nls} with initial data $v(t_0)$ at time $t=t_0$.

We now have to reintroduce the bubble.  As $t_0$ and $t_1$ are far apart, the free evolution of $w(t_0)$ to the future of $t_0$ if $t_1<t_0$
(or to the past of $t_0$ if $t_1>t_0$) is small (of the order $O(\eta_1^{1/100})$).  Choosing $\eta_1$ sufficiently small, an application of
the stability result Lemma~\ref{stab} yields a bound on the scattering size $S(\tilde u)$ either to the future or to the past of $t_0$ and hence,
a bound on either $I_+$ or $I_-$. As Galilean transformations leave the Strichartz norm $S$ invariant and $I_{-}$ and $I_{+}$ were chosen to support
a third of the total $L_{t,x}^{\frac{2(d+2)}{d}}(I_*\times\R^d)$ mass of $u$, choosing $\eta_2$ sufficiently small in \eqref{HUGE}, we derive a
contradiction. Thus, for all $t\in I_0$ there exists $\xi(t)$ such that \eqref{lb-gu} holds.

Using \eqref{lb-gu} (in the same way the the lower bound on the potential energy was used to derive space concentration in the energy-critical
setting), we establish
\begin{proposition}[Space concentration]\label{sc}
Let $u$ be a minimal-mass almost blowup solution on $I_*\times\R^d$ and let $t\in I_0$.  Then, there exist $x(t), \xi(t)\in \R^d$ such that
for any $1<p<\infty$ we have
\begin{align*}
\int_{|x-x(t)|<C(\eta_1)/\tilde N(t)} |G_{\xi(t)}(u)(t)|^p \, dx \geq c(\eta_1) \tilde N(t)^{\frac{pd}2 -d}.
\end{align*}
Here, $\tilde N(t)\in 2^{\Z}$ is the frequency at which the minimal-mass almost blowup solution $G_{\xi(t)}(u)$ is localized
(see Corollary~\ref{freq local}).
\end{proposition}

We omit the details of the proof of Proposition~\ref{sc}.  Given this proposition, in order to prove space localization we need only show

\begin{proposition}\label{sl}
Let $u$ be a minimal-mass almost blowup solution on $I_*\times\R^d$.  Let $t\in I_0$ and let $x(t), \xi(t)\in \R^d$ and $\tilde N(t)\in 2^{\Z}$
be as in Proposition~\ref{sc}.
Then,
\begin{align*}
\int_{|x-x(t)|>\frac{1}{\eta_2\tilde N(t)}} |G_{\xi(t)}(u)(t)|^2 \, dx \leq \eta_1.
\end{align*}
\end{proposition}

The proof of Proposition~\ref{sl} follows very closely that of the analogous statement in the energy-critical case (see \cite{bourg.critical, gopher}).
It is a proof by contradiction.  Using Proposition~\ref{sc}, the conservation of mass, and the pigeonhole principle, one can find a large annulus
where the mass of the solution is small.  One then defines two initial data widely separated in space whose masses are strictly smaller than
the critical mass $m_0$.  One can then find two global solutions to \eqref{nls}, one for each of these two initial data.  In order to use
the stability result Lemma~\ref{stab}, one needs to control the interactions between these two global solutions.
This is done using the pseudoconformal transformation to derive a finite speed of propagation result for the two global solutions.
We omit the details.

\begin{remark}
In the spherically symmetric case, a simple argument based on the conservation of mass shows that $x(t),\xi(t)$ can be taken to be zero
for all $t\in I_0$ (see, for example, \cite{tao} for a similar argument in the energy-critical setting where potential energy is used instead of mass).
\end{remark}

\end{document}